\documentclass[11pt]{article}
\usepackage{t1enc}
\usepackage{lmodern}
\usepackage[T2A,T1]{fontenc}
\usepackage[utf8]{inputenc}
\usepackage[russian,english]{babel}
\usepackage{amsmath,amssymb,amsfonts,amsthm,mathrsfs,textcomp,url,bbm,cancel,comment,enumitem,mathtools}
\usepackage{graphicx}
\usepackage{float}
\usepackage[all,ps]{xy}
\usepackage[colorlinks,urlcolor=cyan,citecolor=MidnightBlue,linkcolor=MidnightBlue]{hyperref}
\usepackage[dvipsnames]{xcolor}
\usepackage[margin=2.5cm]{geometry}
\usepackage{longtable}
\usepackage{tikz}
\makeatletter
\newcommand{\affil}[2]{\gdef\@affila{\textsc{#1}}
  \gdef\@affilb{\textsc{#2}}}
\newcommand{\address}[2]{\gdef\@addressa{#1}
  \gdef\@addressb{#2}}
\newcommand{\email}[2]{\gdef\@emaila{\url{#1}}
  \gdef\@emailb{\url{#2}}}
\newcommand{\@endstuff}{\par\vspace{\baselineskip}\noindent\footnotesize
  \begin{tabular}{@{}l}
    \@affila\\
    \@addressa\\
    \textit{E-mail address:} \@emaila
  \end{tabular}
  \par\medskip\noindent
  \begin{tabular}{@{}l}
    \@affilb\\
    \@addressb\\
    \textit{E-mail address:} \@emailb
  \end{tabular}}
\AtEndDocument{\@endstuff}
\makeatother
\DeclareRobustCommand{\cyrins}[1]{%
  \begingroup\fontfamily{cmr}%
  \foreignlanguage{russian}{#1}%
  \endgroup}
\usetikzlibrary{matrix,calc,decorations,positioning}
\pgfkeys{/tikz/.cd,
  alt double distance/.initial=5pt,
  alt double step/.initial=1pt,}
\pgfdeclaredecoration{double deco}{initial}
{
  \state{initial}[width=\pgfkeysvalueof{/tikz/alt double step},next state=cont] {
    \pgfmoveto{\pgfpoint{\pgfkeysvalueof{/tikz/alt double step}}{\pgfkeysvalueof{/tikz/alt double distance}/2}}
    \pgfpathlineto{\pgfpoint{0.3\pgflinewidth}{\pgfkeysvalueof{/tikz/alt double distance}/2}}
    \pgfpathmoveto{\pgfpoint{0.3\pgflinewidth}{-\pgfkeysvalueof{/tikz/alt double distance}/2}}
    \pgfpathlineto{\pgfpoint{1pt}{-\pgfkeysvalueof{/tikz/alt double distance}/2}}
    \pgfcoordinate{lastup}{\pgfpoint{1pt}{\pgfkeysvalueof{/tikz/alt double distance}/2}}
    \pgfcoordinate{lastdown}{\pgfpoint{1pt}{-\pgfkeysvalueof{/tikz/alt double distance}/2}}}
  \state{cont}[width=\pgfkeysvalueof{/tikz/alt double step}]{
    \pgfmoveto{\pgfpointanchor{lastup}{center}}
    \pgfpathlineto{\pgfpoint{\pgfkeysvalueof{/tikz/alt double step}}{\pgfkeysvalueof{/tikz/alt double distance}/2}}
    \pgfcoordinate{lastup}{\pgfpoint{\pgfkeysvalueof{/tikz/alt double step}}{\pgfkeysvalueof{/tikz/alt double distance}/2}}
    \pgfmoveto{\pgfpointanchor{lastdown}{center}}
    \pgfpathlineto{\pgfpoint{\pgfkeysvalueof{/tikz/alt double step}}{-\pgfkeysvalueof{/tikz/alt double distance}/2}}
    \pgfcoordinate{lastdown}{\pgfpoint{\pgfkeysvalueof{/tikz/alt double step}}{-\pgfkeysvalueof{/tikz/alt double distance}/2}}}
  \state{final}[width=0pt]
  { 
    \pgfmoveto{\pgfpointdecoratedpathlast}}}
\tikzset{alt double/.style={decorate,decoration=double deco}}
\allowdisplaybreaks

\setlist[enumerate]{itemsep=0mm}

\newcommand{\toverset}[2]{%
  \mathop{#2}\limits^{\vbox to -.1ex{\kern-0.4ex\hbox{$\scriptstyle #1$}\vss}}}
\newcommand{\tightoverset}[2]{%
  \mathop{#2}\limits^{\vbox to -.5ex{\kern-0.4ex\hbox{$\scriptstyle #1$}\vss}}}
\renewcommand{\underset}[2]{%
  \mathop{#2}\limits_{\vbox to -.5ex{\kern-1.6ex\hbox{$\scriptstyle #1$}\vss}}}
\newcommand{\tightunderset}[2]{%
  \mathop{#2}\limits_{\vbox to -.5ex{\kern-1.8ex\hbox{$\scriptstyle #1$}\vss}}}
\newcommand{\xra}[1]{%
  \mathop{\xrightarrow{~#1~}}}

\newcommand{\cpt}[1]{%
  \mathop{\toverset{^\bullet}{#1}}}

\newcommand{\usqcup}[1]{%
  \mathop{\tightunderset{#1}{\sqcup}}}
\newcommand{\utimes}[1]{%
  \mathop{\tightunderset{#1}{\times}}}

\def\N{\mathbb{N}}
\def\Z{\mathbb{Z}}
\def\Q{\mathbb{Q}}
\def\R{\mathbb{R}}
\def\C{\mathbb{C}}

\def\RP{\mathbb{R}P}
\def\CP{\mathbb{C}P}

\def\AA{\mathscr{A}}

\def\VV{\mathscr{V}}
\def\SS{\mathscr{S}}

\def\into{\hookrightarrow}
\def\imto{\looparrowright}

\def\ol{\overline}

\def\la{\langle}
\def\ra{\rangle}

\def\im{\mathop{\rm im}}

\def\codim{\mathop{\rm codim}}

\def\id{\mathop{\rm id}}
\def\pr{\mathop{\rm pr}}

\def\Tp{\mathrm{Tp}}

\def\U{\mathrm{U}}

\def\sw1{{w_1{\rm s}}}
\def\sc1{{c_1{\rm s}}}

\def\Emb{\textstyle{\mathop{\rm Emb}}}

\def\Imm{\textstyle{\mathop{\rm Imm}}}

\def\Qhol{\textstyle{\mathop{\rm Qhol}}\displaystyle}

\def\Bihol{\textstyle{\mathop{\rm Bihol}}}
\def\Aut{\textstyle{\mathop{\rm Aut}}}

\def\Sym{\textstyle{\mathop{\rm Sym}}}

{\phantomsection\par\addvspace{.5\baselineskip}\noindent\textbf{Theorem \ref{#1}.\ignorespaces}\begin{em}}%
  {\end{em}\par\addvspace{.5\baselineskip}}%
\newenvironment{prf}[1][\unskip]%
{\par\addvspace{.5\baselineskip}\noindent\textbf{Proof #1.\enspace\ignorespaces}}%
{~$\square$\par\addvspace{.5\baselineskip}}%
\newenvironment{sketchprf}[1][\unskip]%
{\par\addvspace{.5\baselineskip}\noindent\textbf{Sketch of the proof #1.\enspace\ignorespaces}}%
{~$\square$\par\addvspace{.5\baselineskip}}%
{\par\addvspace{.5\baselineskip}\noindent\textit{Claim.\enspace\ignorespaces}\begin{em}}%
  {\end{em}\par\addvspace{.5\baselineskip}}%
{\par\addvspace{.5\baselineskip}\noindent\textit{Proof #1.\enspace\ignorespaces}}%
{~$\diamond$\par\addvspace{.5\baselineskip}}%
{\par\addvspace{.5\baselineskip}\noindent\textit{Remark.\enspace\ignorespaces}}%
{\lpar\addvspace{.5\baselineskip}}%
\newenvironment{subjclass}[1]%
{\par\noindent\begin{small}\textit{#1 Mathematics Subject Classification.\enspace\ignorespaces}}%
  {\end{small}\par}%
\newenvironment{key}%
{\par\noindent\begin{small}\textit{Key words and phrases.\enspace\ignorespaces}}%
  {\end{small}\par}%
{\par\noindent\begin{small}\textbf{Acknowledgement.\enspace\ignorespaces}}%
  {\end{small}\par}%
{\par\addvspace{.5\baselineskip}\noindent\begin{em}}%
  {\end{em}\par\addvspace{.5\baselineskip}}%
\newenvironment{exenum}%
{\begin{ex} ~ \vspace{-.5em}\begin{enumerate}}%
    {\end{enumerate}\end{ex}}%
{\begin{figure}[H]
    \begin{center}
      \centering\includegraphics[scale=0.4]{#1}
      \begin{changemargin}{2cm}{2cm}
        \caption{\footnotesize\hangindent=1.4cm #2}
      \end{changemargin}}%
    {\vspace{-1cm}
    \end{center}
  \end{figure}}%
{\begin{figure}[H]
    \begin{center}
      \centering\includegraphics[scale=0.35]{#1}
      \begin{changemargin}{2cm}{2cm}
        \footnotesize{#2}
      \end{changemargin}}%
    {\vspace{-1cm}
    \end{center}
  \end{figure}}%
\newtheorem{thm}{Theorem}[subsection]%
\newtheorem{lemma}[thm]{Lemma}%
\newtheorem{prop}[thm]{Proposition}%
\newtheorem{crly}[thm]{Corollary}%
\theoremstyle{definition}
\newtheorem{defi}[thm]{Definition}%
\newtheorem{rmk}[thm]{Remark}%
\newtheorem{ex}[thm]{Example}%
\newcommand{\ZeroRoman}[1]{%
\ifcase\value{#1}\relax 0\else\Roman{#1}\fi}

\newcounter{t}

\numberwithin{equation}{section}
\numberwithin{figure}{section}
\title{Quasi-holomorphic maps}
\author{András Csépai \and András Szűcs}
\affil{HUN-REN Alfréd Rényi Institute of Mathematics}{ELTE Eötvös Loránd University, Institute of Mathematics}
\address{Reáltanoda utca 13--15, Budapest, H-1053 Hungary}{Pázmány Péter sétány 1/c, Budapest, H-1117 Hungary}
\email{csepai@renyi.hu}{andras.szucs@ttk.elte.hu}
\date{}
\begin{document}

\maketitle

\begin{abstract}
  We introduce a new notion, called quasi-holomorphic maps. These are real smooth maps equipped with a structure that imitates the singularities and singularity stratifications of holomorphic maps on the source and target manifolds, although the manifolds themselves carry no global complex structures. Some important examples of quasi-holomorphic maps are branched coverings and links of finitely determined holomorphic map germs. We show a Pontryagin--Thom type construction for a ``universal'' quasi-holomor\-phic map with prescribed multisingularities, from which all such maps can be induced, and a similar result for maps with prescribed singularities. Applying this, we prove that the Thom polynomials of holomorphic singularities determine the cohomology classes represented by the singular loci of not only holomorphic but quasi-holomorphic maps as well. As another application we define the cobordism groups of quasi-holomorphic maps with restricted multisingularities, whose classifying space was given by the above construction. We completely compute the free parts of these cobordism groups and in some special cases also obtain results on their torsion parts.
\end{abstract}

\begin{subjclass}{2020}
  57R45 (Primary); 32S20; 
  55Q52; 57R70; 57R90; 58K30 (Secondary)
\end{subjclass}

\begin{key}
  holomorphic map singularities; multisingularities; generalised Pontryagin--Thom construction; singular loci; cobordism groups
\end{key}

\section{Introduction}

The purpose of the present paper is to introduce and investigate a class of real smooth maps equipped with a structure imitating the (multi)singularity stratification of holomorphic maps. The motivation for this is two-fold. Firstly, the notion of \textit{quasi-holomorphic maps} we are about to introduce relaxes the rigidity of holomorphic maps with respect to various deformations without losing global singularity theoretic information. An interesting particular case of this is that cobordism type relations will be well-defined between quasi-holomorphic maps. Secondly, the construction of the ``universal singular map'' (i.e. a map from which any smooth map with specified singularities can be induced) in \cite{rsz} has a complex analogue (i.e. a map from which any holomorphic map with specified singularities can be induced) which was extensively used by Rimányi \cite{rim} to study the global singularity theory of holomorphic maps. It will turn out that quasi-holomorphic maps are exactly the maps classified by this complex universal singular map and they are actually needed for its construction to be completely precise. We note that the precise introduction of quasi-holomorphic maps has been a long-standing open project of the second author; a preliminary version of it appeared in \cite[example 120]{hosszu}.

Throughout the paper by a smooth map we will mean a real $C^\infty$-map between $C^\infty$-manifolds and by a holomorphic map a complex analytic map between complex analytic manifolds. When talking about (co)dimensions, we will always mean real (co)dimensions unless we specify them to be the complex ones. We will sometimes indicate the (real or complex) dimension of a (real or complex) manifold in a superindex but in most cases we will suppress it.

\subsection{Main definitions and examples}
\label{ssec:def}

Before we define quasi-holomorphic maps let us recall a few basic notions in the singularity theory of holomorphic maps. For two holomorphic map germs $\eta,\vartheta\colon(\C^n,0)\to(\C^{n+k},0)$ we say that $\eta$ and $\vartheta$ are \textit{$\AA$-equivalent} if there are biholomorphism germs $\varphi$ and $\psi$ of $(\C^n,0)$ and $(\C^{n+k},0)$ respectively such that $\vartheta=\psi\circ\eta\circ\varphi^{-1}$; further we call the germ $\eta\times\id_\C\colon(\C^{n+1},0)\to(\C^{n+k+1},0)$ the \textit{suspension} of $\eta$. The \textit{singularity} of the germ $\eta$ is the equivalence class of $\eta$ in the equivalence relation generated by $\AA$-equivalence and suspension; it will be denoted by $[\eta]$. The \textit{(complex) codimension} of a singularity $[\eta]$ is the smallest dimension $n$ such that a germ $(\C^n,0)\to(\C^{n+k},0)$ represents $[\eta]$; it will be denoted by $\codim[\eta]$ or shortly $c_\eta$. A \textit{prototype} of a singularity $[\eta]$ is a minimal dimensional germ representing it, i.e. a germ $(\C^{c_\eta},0)\to(\C^{c_\eta+k},0)$ in the equivalence class $[\eta]$. Sometimes, when it is clear from the context, we will abuse notation by omitting the brackets and writing just $\eta$ instead of $[\eta]$.


In the definition of quasi-holomorphic maps we will use the following auxiliary notion.

\begin{defi}
  \label{defi:cstrat}
  A \textit{complex stratification} of a smooth manifold $M^n$ is a collection $\SS=\{(\sigma_i,J_i)\mid i=0,\ldots,\lfloor \frac n2 \rfloor\}$ where $\sigma_i$ is a submanifold of $M$ of codimension $2i$ and $J_i$ is a complex structure on the normal bundle of $\sigma_i$ in $M$ such that the submanifolds $\sigma_i$ are strata of a Whitney stratification on $M$ and the complex structures on the normal bundles of different strata are compatible in the following sense. Observe that for each $i$ and any point $p\in\sigma_i$ the coordinate neighbourhoods of $p$ in $M$ are identified with $\C^i\times\R^{n-2i}$ through $J_i$. We require that for any $j<i$, if $p$ is in $\sigma_i\cap\ol\sigma_j$, then there is a coordinate neighbourhood $U\approx\C^i\times\R^{n-2i}$ of $p$ where $U\cap\sigma_j$ is of the form $V\times\R^{n-2i}$ for a complex submanifold $V\subset\C^i$, moreover, the natural complex structure induced on the normal bundle of $V$ coincides under this identification with the restriction of the complex structure $J_j$.
\end{defi}

\begin{defi}
  \label{defi:qhol}
  A \textit{quasi-holomorphic map} is a triple $(f,\SS,\tilde\SS)$ where $f\colon M^n\to P^{n+2k}$ is a smooth map and $\SS=\{(\sigma_i,J_i)\mid i=0,\ldots,\lfloor \frac n2 \rfloor\}$ and $\tilde\SS=\{(\tilde\sigma_i,\tilde J_i)\mid i=0,\ldots,\lfloor \frac n2 \rfloor+k\}$ are complex stratifications of $M$ and $P$ respectively with the following properties:
  \begin{enumerate}[label=(\roman*)]
  \item for $i=0,\ldots,\lfloor \frac n2 \rfloor$ the restriction $f|_{\sigma_i}$ is a self-transverse immersion into $P$, its image is a dense subset of the closure of $\tilde\sigma_{i+k}$, 
    and for $0<j<k$ the stratum $\tilde\sigma_j$ is empty (note that because of the compatibility of the complex structures $\tilde J_j$ for all $j$, the structure $\tilde J_{i+k}$ naturally extends to the normal bundle of $f|_{\sigma_i}$),
  \item for any point $p\in\sigma_i$ there are coordinate neighbourhoods $U\subset M$ of $p$  and $\tilde U\subset P$ of $f(p)$ identified with $\C^i\times\R^{n-2i}$ and $\C^{i+k}\times\R^{n-2i}$ through $J_i$ and $\tilde J_{i+k}$ respectively such that the germ of $f$ at $p$ in these coordinates is of the form $\eta\times\id_{\R^{n-2i}}$ for a holomorphic map germ $\eta\colon(\C^i,0)\to(\C^{i+k},0)$,
  \item if $p\in\sigma_i$ and $\eta$ are as above, then we have $\codim[\eta]=i$.
  \end{enumerate}
\end{defi}

Although the stratifications $\SS$ and $\tilde\SS$ are part of the quasi-holomorphic structure, we will usually omit them for simplicity and only say that $f\colon M\to P$ is a quasi-holomorphic map.
  
\begin{ex}
  \label{ex:qholmap}
  Naturally, holomorphic maps between complex manifolds are quasi-holomor\-phic.\footnote{For the sake of this example we assume a certain ``genericity'' condition on all maps, which yields that the source and target manifolds are indeed stratified according to the codimensions of the singularities of the maps. We will clarify this assumption in section \ref{sec:sing}.} Other natural examples of quasi-holomorphic maps are the following:
  \begin{enumerate}
  \item Immersions equipped with complex structures on their normal bundles.
  \item Branched coverings: recall that a smooth map $f\colon M^n\to P^n$ is a \textit{branched covering} if there is a submanifold $N^{n-2}\subset M^n$ such that the restriction $f|_{M\setminus N}$ is an $r$-sheeted covering for some $r$ and for each point $p\in N$ there are coordinate neighbourhoods $U\subset M$ of $p$ and $\tilde U\subset P$ of $f(p)$, both identified with $\R^{n-2}\times\C$, where the restriction of $f$ has the form $\id_{\R^{n-2}}\times(-)^r$ (here $(-)^r\colon\C\to\C$ is the map $z\mapsto z^r$).
  \item Restrictions of holomorphic maps $f\colon M\to P$ to preimages of real $C^\infty$-submanifolds $Q\subset P$ that are transverse to all strata determined by the codimensions of the singularities of $f$, i.e. the maps $f|_{f^{-1}(Q)}\colon f^{-1}(Q)\to Q$; in particular the link $S^{2n-1}\to S^{2(n+k)-1}$ of a holomorphic germ $(\C^n,0)\to(\C^{n+k},0)$ is quasi-holomorphic.
  \item Suspensions of holomorphic maps $f\colon M\to P$ to products by real manifolds; that is, if $N$ is a real $C^\infty$-manifold, then the map $f\times\id_N\colon M\times N\to P\times N$ is quasi-holomorphic.
  \end{enumerate}
\end{ex}

\begin{rmk}
  ~ \vspace{-.5em}
  \begin{enumerate}
  \item Quasi-holomorphic maps $f\colon M\to P$ can be slightly perturbed by ambient isotopies, since the complex stratifications of $M$ and $P$ can be carried along the isotopy. This means, for example, that (unlike for holomorphic maps) we can change a quasi-holomorphic map on specific parts, keeping the rest fixed.
  \item The latter two items of example \ref{ex:qholmap} also hold if the map $f$ is only quasi-holomorphic. This already implies that quasi-holomorphic maps can be pulled back, i.e. if $f\colon M\to P$ is a quasi-holomorphic map between smooth manifolds and $g\colon Q\to P$ is a smooth map transverse to all strata $\tilde\sigma_i$ in $P$, then the map $g^*f\colon N\to Q$ indicated on the pullback square
    $$\xymatrix{
      M\ar[r]^{f} & P \\
      N\ar[u]\ar[r]^{g^*f} & Q\ar[u]_g
    }$$
    is quasi-holomorphic as well.
  \end{enumerate}
\end{rmk}


\subsection{Main results and organisation of the paper}
\label{ssec:org}

To precisely state our main results, many notions have to be introduced first. In the following we summarise the content of the paper and intuitively describe the results; references to the precise statements will be written in bold-face font.

Section \ref{sec:sing} consists of preliminary definitions. In subsection \ref{ssec:sing} we define the types of singularities we consider throughout the paper. In subsection \ref{ssec:strat} we describe the stratifications induced by the singularities of holomorphic and quasi-holomorphic maps on the source and target manifolds. Then, in subsection \ref{ssec:taumap} we restrict quasi-holomorphic maps locally so that all of their singularities belong to a given set $\sigma$ of those; such maps are called \textit{$\sigma$-maps}. Sometimes we also restrict maps globally by limiting their possible multiple points in each singular locus; we call this the restriction of their \textit{multisingularities}. Quasi-holomorphic maps with multisingularities restricted to a fixed set $\tau$ of those are called \textit{$\tau$-maps}.

In section \ref{sec:univmap} we recall constructions of ``universal examples'' of holomorphic maps which turn out to have the same universal properties for quasi-holomorphic maps as well. We show in subsections \ref{ssec:globnf} and \ref{ssec:pt} the complex version of the generalised Pontryagin--Thom construction of Rimányi and the second author \cite{rsz} for maps with restricted multisingularities, and clarify the role of quasi-holomorphic maps in the construction. 
The main result here is \textbf{theorem \ref{thm:univmap}} which states that for a set $\tau$ of multisingularities there is a ``universal $\tau$-map'' $f_\tau\colon Y_\tau\to X_\tau$ such that the $\tau$-maps $f\colon M^n\to P^{n+2k}$ are precisely the maps induced from $f_\tau$ through a pullback square
$$\xymatrix{
  Y_\tau\ar[r]^{f_\tau} & X_\tau \\
  M\ar[u]\ar[r]^f & P.\ar[u]
}$$
The existence of such a pullback square was used in \cite{rim} in the case when $f$ is a holomorphic map, however, to precisely construct the map $f_\tau$ the notion of quasi-holomorphic maps is needed. Then in subsection \ref{ssec:kaz}, corresponding to a set $\sigma$ of singularities without global restrictions, we show the construction of a strongly related space $K_\sigma$, called 
the \textit{Kazarian space} (see e.g. \cite{kazspace}). We show in \textbf{theorem \ref{thm:kaz}} that the singular strata of a $\sigma$-map $f\colon M\to P$ are transverse preimages of specified subspaces of the Kazarian space under a map $M\to K_\sigma$; this was known for holomorphic maps, we will see that it holds for quasi-holomorphic ones as well.

Section \ref{sec:globsingth} introduces two of the main areas of global singularity theory specified to quasi-holomorphic maps, applying the previous universal constructions. 
First, in subsection \ref{ssec:loci} we observe that the spaces $Y_\tau$ and $K_\sigma$ correspond to the study of (multi)singular loci of quasi-holomorphic maps. This immediately yields \textbf{theorem \ref{thm:tp}}, stating that the Poincaré dual cohomology classes represented by the singular loci of quasi-holomorphic maps are determined by the \textit{Thom polynomials} of the singularities in exactly the same way as for holomorphic maps. Then, in subsection \ref{ssec:cob} we turn to the space $X_\tau$; we define the \textit{cobordism group} of $\tau$-maps from $n$-manifolds to a fixed manifold $P^{n+2k}$, denoted by $\Qhol_\tau(P)$, and show in \textbf{theorem \ref{thm:cob}} that we have $\Qhol_\tau(P)\cong[\cpt P,X_\tau]$ where $\cpt P$ denotes the one-point compactification of $P$. That is, $X_\tau$ is the \textit{classifying space} of cobordisms of $\tau$-maps.

Then in section \ref{sec:cob} we apply theorem \ref{thm:cob} to study some concrete groups $\Qhol_\tau(P)$. The main tool here is introduced in subsection \ref{ssec:key}; it is the complex version of the \textit{key fibration} in \cite{hosszu} which connects classifying spaces of cobordisms of maps with different singularity restrictions. The existence of the real smooth version of this fibration has a proof that works in the quasi-holomorphic setting without change; see \textbf{theorem \ref{thm:key}}. By this we have that if $\sigma$ is a set of singularities and $\eta\in\sigma$ is the most complicated 
element, further $\tau$ and $\tau'$ are the sets of all multisingularities composed of the elements of $\sigma$ and $\sigma\setminus\{\eta\}$ respectively, then there is a fibration
$$X_\tau\xra{X_{\tau'}}\Gamma T\tilde\xi_\eta$$
where $\tilde\xi_\eta$ is a vector bundle appearing in the construction of the universal $\tau$-map $f_\tau$, $T\tilde\xi_\eta$ is its Thom space and $\Gamma$ is the infinite loop space of infinite suspension functor $\Omega^\infty S^\infty$. We note that for a vector bundle $\zeta$, the space $\Gamma T\zeta$ is the classifying space of the cobordism groups $\Imm^\zeta(P)$ of immersions into some fixed target $P$ with normal bundles induced from $\zeta$.

In subsection \ref{ssec:rat} we show that the key fibration is rationally trivial for quasi-holomorphic $\sigma$-maps, yielding \textbf{theorem \ref{thm:keyractrivi}} which states that the rationalisation of $X_\tau$ is a product of 
classifying spaces of cobordisms of immersions with normal bundles corresponding to the singularities in $\sigma$. Its main use is \textbf{corollary \ref{crly:ractrivi}} stating that for any manifold $P$ the vector space $\Qhol_\tau(P)\otimes\Q$ is the direct sum of the vector spaces $\Imm^{\tilde\xi_\eta}(P)\otimes\Q$ for all $\eta\in\sigma$. Geometrically this means that up to rational cobordism quasi-holomorphic maps behave as if their singular strata were independently immersed. Then in subsection \ref{ssec:fold} we obtain information on the torsion parts of cobordism groups of quasi-holomorphic \textit{fold} maps of $n$-manifolds to $\R^{n+2}$, denoted by $\Qhol_{\mu_1}(\R^{n+2})=:\Qhol_{\mu_1}(n,1)$. 
Analysing the Thom space of the vector bundle $\tilde\xi_\eta$, where $\eta$ is the fold singularity, yields \textbf{theorem \ref{thm:foldkey}} by which there is a long exact sequence
$$\ldots\to\pi^s_{n+2}(\CP^\infty)\to\Qhol_{\mu_1}(n,1)\to\pi^s_{n+2}(\CP^\infty/\RP^\infty\wedge\CP^\infty/\CP^1)\to\pi^s_{n+1}(\CP^\infty)\to\ldots$$
This implies \textbf{corollary \ref{crly:foldtors}} on the vanishing of some of the torsion part of the group $\Qhol_{\mu_1}(n,1)$; namely this torsion group has no prime components for primes greater than $\frac{n+5}2$.

Finally, in section \ref{sec:fin} we pose open problems and possible further directions for future works in the investigation of quasi-holomorphic maps and their applications.

\section{Singularities and singular maps}
\label{sec:sing}


From now on let an integer $k\ge0$ be fixed. In this section we will describe singularity structures of quasi-holomorphic maps of fixed codimension $2k$. This will mainly consist of classical definitions and examples applied in the quasi-holomorphic setting.

\subsection{Stable and finitely determined singularities}
\label{ssec:sing}

In the following we define singularities that are ``nice'' in some appropriate sense.

\begin{defi}
  For two complex analytic manifolds $M$ and $P$ the product of the biholomorphism groups of $M$ and $P$ acts on the space of proper holomorphic maps $f\colon M\to P$ defined by the formula $(\varphi,\psi)\mapsto\psi\circ f\circ\varphi^{-1}$. The map $f$ is said to be \textit{stable} if it is in the interior of an orbit of this action, that is, if it has a neighbourhood in which every holomorphic map $g$ is such that there are biholomorphisms $\varphi$ of $M$ and $\psi$ of $P$ for which $g$ is $\psi\circ f\circ\varphi^{-1}$.
\end{defi}

\begin{defi}
  Let $\eta\colon(\C^n,0)\to(\C^{n+k},0)$ be a holomorphic map germ.
  \begin{enumerate}
  \item We call $\eta$ a \textit{stable} germ if it is the germ of a stable map around the origin.
  \item We call $\eta$ a \textit{finitely determined} germ if there are neighbourhoods $U$ and $\tilde U$ of the origin in $\C^n$ and $\C^{n+k}$ respectively and a representative $f\colon U\to\tilde U$ of $\eta$ such that the restriction $f|_{U\setminus\{0\}}\colon U\setminus\{0\}\to\tilde U\setminus\{0\}$ is a stable map.
  \end{enumerate}
\end{defi}

\begin{rmk}
  The above definitions of stable and finitely determined germs are not standard, however, they are equivalent to the standard ones. For this notion of stability see the infinitesimal criterion for proper maps in \cite{math5} and the preceding papers of Mather in this series, and for finite determinacy see the Mather--Gaffney criterion, e.g. \cite[theorem 2.1]{wallfin}.
\end{rmk}

Since the stability and finite determinacy of germs are independent of the coordinates, the following definition is correct.

\begin{defi}
  A singularity $[\eta]$ of holomorphic germs of codimension $k$ is said to be \textit{stable} (resp. \textit{finitely determined}) if its prototype is a stable (resp. finitely determined) germ $(\C^{c_\eta},0)\to(\C^{c_\eta+k},0)$.
\end{defi}

\begin{exenum}
  \item The ``non-singular singularity'', i.e. the singularity class $\Sigma^0:=[(\C^0,0)\into(\C^k,0)]$ of the regular germs is clearly stable.
  \item For any $r\ge1$ the singularity class of the $r$'th power map $(-)^r\colon(\C,0)\to(\C,0)$ is finitely determined. However, it can be shown that it is not stable if $r\ge3$. The stable singularity $[(-)^2]$ is called the ($0$-codimensional) \textit{fold}.
\end{exenum}

Throughout this paper we only consider maps every singularity of which is finitely determined; we will tacitly assume this from now on. This assumption was also present in example \ref{ex:qholmap}; it was the ``genericity'' condition mentioned in the footnote.

\begin{rmk}
  Unfortunately, despite the above assumption, this finite determinacy condition is not always true for \textit{almost all maps} in the sense that the space of maps satisfying it is not always dense in the space of all maps. 
  However, there is a range of dimensions in which this finite determinacy 
  does hold for almost all maps; see e.g. \cite[theorem 5.6]{wallfin}. Moreover, in another (smaller) range of dimensions it even holds that almost all maps are stable; see e.g. \cite[theorem 5.5]{wallfin}. We also note that ``small perturbations'' of stable maps have the same singularities as the original maps; this is not true for all maps in general.
\end{rmk}

\subsection{Singular strata}
\label{ssec:strat}

\begin{defi}
  For a singularity $[\eta]$ of finitely determined holomorphic germs of complex codimension $k$ and a quasi-holomorphic map $f\colon M^n\to P^{n+2k}$ we call a point $p\in M$ an \textit{$\eta$-point} if the germ of $f$ at $p$ belongs to $\eta$ in some (and hence any) local coordinate systems. The set of $\eta$-points in $M$ is denoted by $\eta(f)$ and called the \textit{$\eta$-locus} of $f$; the image of the $\eta$-locus in $P$ is denoted by $\tilde\eta(f):=f(\eta(f))$.
\end{defi}

It is not hard to see that in the above setting the locus $\eta(f)\subset M$ is a submanifold of codimension $2c_\eta$ and the restriction $f|_{\eta(f)}\colon\eta(f)\imto P$ is an immersion. Moreover, we also have that the submanifolds $\eta(f)$ stratify the manifold $M$ and this stratification refines the codimension stratification formed by the submanifolds $\sigma_i\subset M$ in the quasi-holomorphic structure of the map $f$.

In the rest of this subsection we define a few examples to be used later. An important classification of singularities is by the corank of their differentials, as follows. For any $r\ge0$, we define $\Sigma^r$ to be the set of all singularities $[\eta\colon(\C^n,0)\to(\C^{n+k},0)]$ for which we have $\dim_\C\ker d\eta_0=r$. For a quasi-holomorphic map $f\colon M^n\to P^{n+2k}$ we denote by $\Sigma^r(f)$ the union of the sets $\eta(f)$ for all $\eta\in\Sigma^r$. Then for almost all maps $f$ the set $\Sigma^r(f)$ is a submanifold in $M$ of codimension $2r(k+r)$ (see e.g. \cite{boar}). If $r\ge2$, then the structures of the singularities in $\Sigma^r$ are complicated and not understood in general. However, at least the \textit{stable} singularities in $\Sigma^1$ are well-known: they were classified by Morin \cite{mor}, who gave a canonical form of their prototypes; these singularities belong to a single sequence of singularities $A_i$ of (complex) codimension $i(k+1)$, parametrised by the positive integers. The same paper of Morin also proves that for almost all maps $f$ there are no unstable singularities of $f$ in $\Sigma^1$.

\begin{defi}
  The stable singularity $A_i$ 
  is called the $i$'th \textit{Morin singularity}. A quasi-holomor\-phic map whose every singularity is a Morin singularity (i.e. a stable map $f$ with $\Sigma^r(f)=\varnothing$ for $r\ge2$) is called a quasi-holomorphic \textit{Morin map}.
\end{defi}

\begin{exenum}
\item The singularity $A_1$ (of $k$-codimensional germs) is called the ($k$-codimensional) \textit{fold}. Its prototype is
  $$\C^{k+1}\to\C^{2k+1};~(x_1,\ldots,x_k,y)\mapsto(x_1,\ldots,x_k,x_1y,\ldots,x_ky,y^2).$$
\item The singularity $A_2$ (of $k$-codimensional germs) is called the ($k$-codimensional) \textit{cusp}. Its prototype is
  $$\C^{2k+2}\to\C^{3k+2};~(x_1,\ldots,x_{2k+1},y)\mapsto(x_1,\ldots,x_k,x_1y+x_2y^2,\ldots,x_{2k-1}y+x_{2k}y^2,x_{2k+1}y+y^3).$$
\item In the $k=0$ case, fold is the singularity of the squaring map $(-)^2$. For $r\ge 3$ the $r$'th power singularity $[(-)^r]$ also belongs to $\Sigma^1$ but it is not stable.
\end{exenum}

\subsection{Maps with restricted multisingularities}
\label{ssec:taumap}

So far, when talking about singularities, we only considered the local structure of maps. We now introduce a global property.

\begin{defi}
  A \textit{multisingularity} is a formal finite sum $\mu=m_1[\eta_1]+\ldots+m_r[\eta_r]$ of singularities of finitely determined holomorphic germs $\eta_i\colon(\C^n,0)\to(\C^{n+k},0)$. We say that a multigerm $(\C^n,S)\to(\C^{n+k},0)$ (where $S\subset\C^n$ is a finite set) \textit{belongs to the multisingularity} $\mu$ if the set $S$ consists of $m_1+\ldots+m_r$ points out of which precisely $m_i$ are $\eta_i$-points for each $i$.
\end{defi}

\begin{defi}
  For a quasi-holomorphic map $f\colon M\to P$ and a multisingularity $\mu$ we call a point $p\in M$ a \textit{$\mu$-point} if the multigerm of $f$ at $f^{-1}(f(p))$ belongs to $\mu$ in some (and hence any) local coordinate systems. The set of $\mu$-points in $M$ is denoted by $\mu(f)$ and called the \textit{$\mu$-locus} of $f$; the image of the $\mu$-locus in $P$ is denoted by $\tilde\mu(f):=f(\mu(f))$.
\end{defi}

\begin{rmk}
  For a multisingularity $\mu=m_1[\eta_1]+\ldots+m_r[\eta_r]$ and a quasi-holomorphic map $f\colon M\to P$ the restriction $f|_{\mu(f)}\colon\mu(f)\to\tilde\mu(f)$ is an $(m_1+\ldots+m_r)$-fold covering. Also note that $\tilde\mu(f)$ is the intersection for all $i$ of the $m_i$-fold self-intersection of $\tilde\eta_i(f)$.
\end{rmk}



In the following we put both local and global restrictions on quasi-holomorphic maps by restricting their possible multisingularities.


\begin{defi}
  Let $\tau$ be a set of multisingularities. We say that a quasi-holomorphic map $f\colon M^n\to P^{n+2k}$ is a \textit{$\tau$-map} if each multigerm of $f$ belongs to a multisingularity in $\tau$.
\end{defi}

\begin{exenum}
\item $\{1\Sigma^0\}$-maps are embeddings with complex structures on their normal bundles.
\item More generally, $\{i\Sigma^0\mid i=1,\ldots,r\}$-maps for some $r$ are immersions with at most $r$-tuple points, again with complex normal bundles.
\item $\{1A_1\}\cup\{i\Sigma^0\mid i=1,2,\ldots\}$-maps are quasi-holomorphic \textit{simple fold maps} where ``simple'' means that singular points are not multiple. 
\item In the $k=0$ case, $\{r\Sigma^0,1[(-)^r]\}$-maps are the $r$-sheeted \textit{simple branched covernigs}.
\end{exenum}

\begin{defi}
  Let $\tau$ be a set of multisingularities. If $\tau$ consists of all linear combinations of the elements of a set of singularities $\sigma$ (that is, $\tau=\N\la\sigma\ra$), then we call $\tau$ a \textit{complete} set of multisingularities. In this case we call $\tau$-maps also $\sigma$-maps, that is, a \textit{$\sigma$-map} is a quasi-holomorphic map, each singularity of which belongs to $\sigma$.
\end{defi}

\begin{exenum}
\item $\{\Sigma^0\}$-maps are immersions with complex normal bundles.
\item $\{\Sigma^0,A_1\}$-maps and $\{\Sigma^0,A_1,A_2\}$-maps are the quasi-holomorphic \textit{fold maps} and \textit{cusp maps} respectively.
\item Stable $\Sigma^0\cup\Sigma^1$-maps, i.e. $\Sigma^0\cup\{A_i\mid i=1,2,\ldots\}$-maps are the quasi-holomorphic Morin maps.
\item In the $k=0$ case, $\{\Sigma^0,[(-)^r]\}$-maps are the $r$-sheeted branched coverings.
\end{exenum}


There is a natural partial order of multisingularities (resp. singularities) defined by $\mu\ge\nu$ (resp. $[\eta]\ge[\vartheta]$) iff in any neighbourhood of a $\mu$-point (resp. $\eta$-point) there is a $\nu$-point (resp. $\vartheta$-point). In this case we also say that $\mu$ (resp. $[\eta]$) is \textit{more complicated} than $\nu$ (resp. $[\vartheta]$) and $\nu$ (resp. $[\vartheta]$) is \textit{simpler} than $\mu$ (resp. $[\eta]$). Clearly if $f$ is a $\tau$-map (resp. $\sigma$-map), then $f$ is also a $\tau'$-map (resp. $\sigma'$-map) where $\tau'$ (resp. $\sigma'$) is the maximal downward closed subset of $\tau$ (resp. $\sigma$), hence from now on we will always tacitly assume that (multi)singularity sets are downward closed.

\begin{exenum}
\item For regular multigerms, if $k\ge1$, we have $1\Sigma^0<2\Sigma^0<\ldots$
\item Morin singularities form a single line $A_1<A_2<\ldots$
\item A single fold point is more complicated than a double regular point, i.e. $1A_1>2\Sigma^0$.
\end{exenum}


\section{Universal constructions 
}
\label{sec:univmap}

In this section we recall the complex version of the extended Pontryagin--Thom construction in \cite{rsz}; see also \cite{univmap} for an earlier version containing the main ideas, \cite{rim} for applications of the complex version and \cite{rlsym} for algebraic geometric preliminaries. This construction yields, for a set $\tau$ of multisingularities of finitely determined holomorphic germs of complex codimension $k\ge0$, the 
``universal $\tau$-map'' $f_\tau\colon Y_\tau\to X_\tau$. Its main ingredient is the definition of a \textit{global normal form} for each multisingularity $\mu\in\tau$, by which we mean a universal example for $\tau$-maps in a neighbourhood of the $\mu$-locus; then the universal $\tau$-map is ``glued'' together from these global normal forms. Another way of looking at this, which was utilised by Kazarian, is that we ``cut'' the universal map into the global normal forms; see e.g. \cite{kaztp}. His approach results in another ``universal'' space $K_\tau$, which is particularly useful if $\tau$ is a complete multisingularity set and whose construction does not involve gluings, hence some properties of it are easier to understand.

Throughout this section we shall not give complete proofs of the results except for the ones necessary to place quasi-holomorphic maps into the construction. 

\begin{rmk}
  In \cite{rsz} the Pontryagin--Thom type construction is only stated if $\tau$ is a set of multisingularities of \textit{stable} germs, since this was the case they wanted to apply it to. However, this stability was only needed when applying \cite[lemma 3]{rsz} which also holds for the larger class of \textit{finitely determined} germs, hence the construction holds in this case as well.
\end{rmk}

\subsection{Global normal forms}
\label{ssec:globnf}

For more details on the following constructions see works of Rimányi, e.g. \cite{rim}. Suppose that $\eta\colon(\C^{c_\eta},0)\to(\C^{c_\eta+k},0)$ is a finitely determined holomorphic germ which is a prototype of its singularity class $[\eta]$. The symmetry group
$$\Aut_\AA\eta:=\{(\varphi,\psi)\in\Bihol_0(\C^{c_\eta})\times\Bihol_0(\C^{c_\eta+k})\mid\psi\circ\eta\circ\varphi^{-1}=\eta\}$$
of the germ $\eta$ has a maximal compact subgroup $G_\eta$. For a multisingularity $\mu=m_1[\eta_1]+\ldots+m_r[\eta_r]$ these groups corresponding to the prototypes $\eta_i$ of the singularities in $\mu$ form the maximal compact subgroup $G_\mu=\toverset{r}{\underset{i=1}\prod}(\Sym_{m_i}\ltimes G_{\eta_i}^{m_i})$ of the symmetry group of $\mu$. The representations of this group on the source and target spaces of multigerms belonging to $\mu$ yield vector bundles $\xi_\mu\to B_\mu$ and $\tilde\xi_\mu\to\tilde B_\mu$ respectively together with an $m$-sheeted covering map $\varphi_\mu\colon B_\mu\to\tilde B_\mu$, where $\tilde B_\mu$ is the classifying space $BG_\mu$ and $m=m_1+\ldots+m_r$. Taking in all fibres of $\xi_\mu$ the multigerms representing $\mu$ gives a map $\Phi_\mu\colon\xi_\mu\to\tilde\xi_\mu$ which is a fibrewise (non-linear) map of vector bundles
$$\xymatrix{
  \xi_\mu\ar[r]^{\Phi_\mu}\ar[d] & \tilde\xi_\mu\ar[d] \\
  B_\mu\ar[r]^{\varphi_\mu} & \tilde B_\mu.
}$$
The universal $\tau$-map $f_\tau\colon Y_\tau\to X_\tau$ will be built out of the blocks $\Phi_\mu|_{D\xi_\mu}\colon D\xi_\mu\to D\tilde\xi_\mu$ corresponding to the elements $\mu$ of $\tau$ (where $D\zeta$ denotes the disk bundle of the vector bundle $\zeta$).


\begin{lemma}
  \label{lemma:semig}
  For any quasi-holomorphic map $f\colon M\to P$ and any multisingularity $\mu$, if $f$ has no more complicated multisingularities than $\mu$, then there are tubular neighbourhoods $U_\mu\subset M$ and $\tilde U_\mu\subset P$ of $\mu(f)$ and $\tilde\mu(f)$ 
  respectively that fit into a commutative diagram
  $$\xymatrix@R=.333pc{   
    \xi_\mu\ar[rrr]^{\Phi_\mu}\ar[ddddd] &&& \tilde\xi_\mu\ar[ddddd] \\
    & U_\mu\ar[ddd]\ar[r]^{f|_{U_\mu}}\ar[ul] & \tilde U_\mu\ar[ddd]\ar[ur] &\\ \\ \\
    & \mu(f)\ar[r]^{f|_{\mu(f)}}\ar[dl] & \tilde\mu(f)\ar[dr] &\\
    B_\mu\ar[rrr]^{\varphi_\mu} &&& \tilde B_\mu
  }$$
  where the left-hand and right-hand inner squares are pullback squares (the vertical arrows are the projections).
\end{lemma}

\begin{prf}
  We will only prove this if $\mu$ consists of a single monosingularity $[\eta]$ (i.e. $\mu=1[\eta]$); for general multisingularities the proof is analogous, only formally much more complicated. Now $f|_{\mu(f)}\colon\mu(f)\to\tilde\mu(f)$ is a diffeomorphism (for a general multisingularity it would be an $m$-sheeted covering where $m=m_1+\ldots+m_r$ is the multiplicity in $\mu$), hence we can put $N:=\mu(f)=\tilde\mu(f)$ as abstract manifolds. We will suppose that we have $\mu(f)\subset\sigma_i$ (i.e. $\codim[\eta]=i$) and we choose tubular neighbourhoods $U_\mu\subset M\setminus\{\sigma_j\mid j>i\}$ and $\tilde U_\mu\subset P\setminus\{\tilde\sigma_{j+k}\mid j>i\}$ (which exist because $\mu(f)$ and $\tilde\mu(f)$ are proper); see definition \ref{defi:qhol}.
  
  Now the germs of $U_\mu$ and $\tilde U_\mu$ at $N$ are bundle germs with structure groups $\pr(\Aut_\AA\eta)$ and $\tilde\pr(\Aut_\AA\eta)$ respectively where $\pr$ and $\tilde\pr$ denote the projections of the group $\Aut_\AA\eta\subset\Bihol_0(\C^i)\times\Bihol_0(\C^{i+k})$ to the first and the second factor respectively; by this we mean that there is an open cover $\VV$ of $N$ such that for each $V\in\VV$ we have $U_\mu|_V\approx V\times(\C^i,0)$ and $\tilde U_\mu|_V\approx V\times(\C^{i+k},0)$ and the transition maps are projected images of elements of $\Aut_\AA\eta$. Moreover, these bundle germs are associated to each other in the sense of \cite[lemma 4]{rsz} (actually, the complex version of this lemma), that is, for any $V,V'\in\VV$ the transition maps of $U_\mu$ and $\tilde U_\mu$ at $V\cap V'$ are of the form $\pr\circ\varphi_{V,V'}$ and $\tilde\pr\circ\varphi_{V,V'}$ for a ``holomorphic'' map $\varphi_{V,V'}\colon V\cap V'\to\Aut_\AA\eta$ respectively; here ``holomorphic'' means that $V\cap V'$ can be covered by open subsets $W$ such that the maps $\id_W\times(\pr\circ\varphi_{V,V'}|_W)(-)$ and $\id_W\times(\tilde\pr\circ\varphi_{V,V'}|_W)(-)$ are biholomorphism germs at $W$ of $W\times(\C^i,0)$ and $W\times(\C^{i+k},0)$ respectively.

  Now by \cite[lemma 4]{rsz} the structure groups of $U_\mu$ and $\tilde U_\mu$ can be reduced to the projections of the maximal compact subgroup $G_\eta$ of $\Aut_\AA\eta$ such that they remain associated to each other. Hence $U_\mu$ and $\tilde U_\mu$, which can be identified with the normal bundles of $\mu(f)$ and $\tilde\mu(f)$, are pullbacks of $\xi_\eta\to BG_\eta=B_\mu$ and $\tilde\xi_\eta\to BG_\eta=\tilde B_\mu$ respectively and these pullbacks fit into a commutative diagram as claimed.
\end{prf}

\begin{rmk}
  If $f\colon M\to P$ is a general quasi-holomorphic map, then we can iterate lemma \ref{lemma:semig} with respect to the multisingularities of $f$, descending from the higher to the lower ones. At each step, corresponding to some multisingularity $\mu$, this produces tubular neighbourhoods $U_\mu$ and $\tilde U_\mu$ in the complement of the previously constructed neighbourhoods of the higher codimensional strata. This procedure partitions $M$ and $P$ to tubular neighbourhoods of multisingular loci with given vector bundle structures.
\end{rmk}

Lemma \ref{lemma:semig} has the following converse statement:

\begin{lemma}
  \label{lemma:pb}
  Let $\mu$ be a multisingularity and $K$, $L$ manifolds and $U$, $V$ total spaces of vector bundles over them that fit into a diagram
  $$\xymatrix@R=.333pc{   
    \xi_\mu\ar[rrr]^{\Phi_\mu}\ar[ddddd] &&& \tilde\xi_\mu\ar[ddddd] \\
    & U\ar[ddd]\ar[r]^f\ar[ul] & V\ar[ddd]\ar[ur] &\\ \\ \\
    & K\ar[r]^{f|_K}\ar[dl] & L\ar[dr] &\\
    B_\mu\ar[rrr]^{\varphi_\mu} &&& \tilde B_\mu
  }$$
  where the left-hand, right-hand and bottom inner squares are pullback diagrams. Then $f$ has a natural quasi-holomorphic map structure that has multisingularity $\mu$ at $K$ and below $\mu$ everywhere else.
\end{lemma}

\begin{prf}
  The map $\Phi_\mu$ restricted to the preimage of the image of each fibre of $\xi_\mu$ has multigerm at zero in the multisigularity $\mu$, hence $\xi_\mu$ and $\tilde\xi_\mu$ are stratified according to the multisingularities at most $\mu$. The preimages of these strata form stratifications of $U$ and $V$ which have complex normal bundles induced from those in $\xi_\mu$ and $\tilde\xi_\mu$ and the map $f$ with this structure is clearly quasi-holomorphic.
\end{prf}

The above two lemmas show that the map $\Phi_\mu\colon\xi_\mu\to\tilde\xi_\mu$ is 
the canonical universal example of the mapping of a tubular neighbourhood of the locus of $\mu$-points in the source to a tubular neighbourhood of its image for any quasi-holomorphic map. This justifies to call $\Phi_\mu$ the \textit{global normal form} of the multisingularity $\mu$.

\subsection{Extended Pontryagin--Thom construction}
\label{ssec:pt}

By piecing together the global normal forms of all elements in a multisingularity set $\tau$ we obtain the universal example for quasi-holomorphic $\tau$-maps, as the following theorem shows. This is a direct analogue of \cite[theorem 1]{rsz} and its proof is completely analogous to it, hence we will only sketch its construction, by which we also define the spaces $X_\tau$, $Y_\tau$ and the map $f_\tau$.

\begin{thm}
  \label{thm:univmap}
  For any set $\tau$ of multisingularities of finitely determined holomorphic multigerms, there is a map $f_\tau\colon Y_\tau\to X_\tau$ with the following properties:
  \begin{enumerate}
  \item For any quasi-holomorphic $\tau$-map $f\colon N^{n-1}\to\partial P^{n+2k}$ from a closed manifold $N^{n-1}$ to the boundary of a manifold $P^{n+2k}$, if there are maps $\kappa_N$ and $\kappa_P$ such that the outer square in the diagram below is a pullback diagram, then there is a manifold $M^n$ with boundary $\partial M=N$, an extension $\kappa_M$ of $\kappa_N$ and a quasi-holomorphic $\tau$-map $F\colon M\to P$ such that the upper inner square is a pullback diagram as well
    $$\xymatrix{
      Y_\tau\ar[rr]^{f_\tau} && X_\tau\\
      & M\ar@{-->}[ul]_(.35){\kappa_M}\ar@{-->}[r]^F & P\ar[u]_{\kappa_P}\\
      N\ar[uu]^{\kappa_N}\ar@{^(-->}[ur]\ar[rr]^f && \partial P.\ar@{^(->}[u]
    }$$
  \item For any quasi-holomorphic $\tau$-map $F\colon M^n\to P^{n+2k}$ between manifolds with boundaries, if there are maps $\kappa_{\partial M}$ and $\kappa_{\partial P}$ such that the outer square in the diagram below is a pullback diagram, then there are extensions $\kappa_M$ and $\kappa_P$ of $\kappa_{\partial M}$ and $\kappa_{\partial P}$ respectively such that the upper inner square is a pullback diagram as well
    $$\xymatrix{
      Y_\tau\ar[rrr]^{f_\tau} &&& X_\tau\\
      & M\ar@{-->}[ul]_(.35){\kappa_M}\ar[r]^F & P\ar@{-->}[ur]^(.35){\kappa_P} &\\
      \partial M\ar[uu]^{\kappa_{\partial M}}\ar@{^(->}[ur]\ar[rrr]^{F|_{\partial M}} &&& \partial P\ar@{_(->}[ul]\ar[uu]_{\kappa_{\partial P}}.
    }$$
  \end{enumerate}
\end{thm}

In particular, if the manifold $P$ above has empty boundary, then this theorem states that $\tau$-maps $F\colon M\to P$ correspond to pullback squares inducing $F$ from $f_\tau$. This shows that the map $f_\tau$ deserves to be called the ``universal $\tau$-map''.

\begin{sketchprf}
  We proceed by (possibly transfinite) induction on an arbitrary complete order extending the natural partial order of the multisingularities in $\tau$. The base case is $\tau=\{1\Sigma^0\}$, i.e. the case of embeddings with complex structures on their normal bundles. In this case the usual Pontryagin--Thom construction states that the embedding $B\U(k)\into T\gamma_k$ of the zero-section into the Thom space of the tautological complex $k$-vector bundle serves as the map $f_{\{1\Sigma^0\}}$.
  
  For the induction step let $\mu\in\tau$ be the maximal element and put $\tau':=\tau\setminus\{\mu\}$. Restricting the map $\Phi_\mu$ to the sphere bundle of $\xi_\mu$ is a map $\Phi_\mu|_{S\xi_\mu}\colon S\xi_\mu\to S\tilde\xi_\mu$ whose every multisingularity is below $\mu$, in particular it is a $\tau'$-map. Hence by the induction hypothesis there are maps $\rho_\mu\colon S\xi_\mu\to Y_{\tau'}$ and $\tilde\rho_\mu\colon S\tilde\xi_\mu\to X_{\tau'}$ inducing it from $f_{\tau'}$. We put
  $$Y_\tau:=Y_{\tau'}\usqcup{\rho_\mu}D\xi_\mu,\quad X_\tau:=X_{\tau'}\usqcup{\tilde\rho_\mu}D\tilde\xi_\mu\quad\text{and}\quad f_\tau:=f_{\tau'}\cup\Phi_\mu|_{D\xi_\mu}$$
  which completes the recursive construction of the map $f_\tau$. To prove the two claims we use the universal property of the map $f_{\tau'}$ guaranteed by the induction hypothesis and that of the map $\Phi_\mu$ guaranteed by lemmas \ref{lemma:semig} and \ref{lemma:pb} and the fact that $f_\tau$ was constructed using the same universal property of $f_{\tau'}$.
\end{sketchprf}

\begin{rmk}
  In the construction of the map $f_\tau$ sketched above we essentially use the notion of quasi-holomorphicity as opposed to holomorphicity, since the fibres of the sphere bundles $S\xi_\mu$ and $S\tilde\xi_\mu$ are odd dimensional submanifolds in the corresponding fibres of $\xi_\mu$ and $\tilde\xi_\mu$, hence the restriction $\Phi_\mu|_{S\xi_\mu}$ is only quasi-holomorphic, not holomorphic.
\end{rmk}

\subsection{Kazarian's space}
\label{ssec:kaz}

The following constructions are explained in detail in e.g. \cite{kaztp} and \cite{hosszu}. In this subsection let $\sigma$ be a set of singularities without multiplicities. Then one can consider in the infinite jet space $J_0^\infty(\C^n,\C^{n+k})$ the subspace $J_\sigma(n)$ consisting of the jets of germs which belong to singularities in $\sigma$. We will denote by $\gamma_n\to B\U(n)$ the tautological complex $n$-vector bundle and by $J_0^\infty(\gamma_n,\gamma_{n+k})$ the infinite jet bundle over $B\U(n)\times B\U(n+k)$ of bundle maps $\gamma_n\to\gamma_{n+k}$ (which has $J_0^\infty(\C^n,\C^{n+k})$ as fibre).

\begin{defi}
  \label{defi:kaz}
  For a set $\sigma$ of singularities of finitely determined holomorphic germs, the \textit{Kazarian space} of $\sigma$-maps of $n$-manifolds to $(n+k)$-manifolds is the total space $\hat K_\sigma(n)$ of the subbundle of $J_0^\infty(\gamma_n,\gamma_{n+k})\to B\U(n)\times B\U(n+k)$ whose fibre is $J_\sigma(n)$.
\end{defi}

The Kazarian space is naturally the union of its subspaces corresponding to the individual singularities in $\sigma$. More concretely, for an element $[\eta]\in\sigma$, the union in all fibres of $J_0^\infty(\gamma_n,\gamma_{n+k})$ the subspace of $J_0^\infty(\C^n,\C^{n+k})$ consisting of jets of germs that belong to $[\eta]$, forms a subset $B_\eta(n)\subset\hat K_\sigma(n)$ whose closure is a subvariety. This is connected to the global normal forms of singularities in that we have $B_\eta(n)\cong BG_\eta$ (provided $n\ge c_\eta$), moreover, the normal bundle of $B_\eta(n)$ in $\hat K_\sigma(n)$ is $\xi_\eta 
$. The main use of the Kazarian space space is that if $f\colon M^n\to P^{n+k}$ is a \textit{holomorphic} $\sigma$-map, then there is a map $\hat\kappa_f\colon M\to\hat K_\sigma(n)$ such that for all singularities $[\eta]\in\sigma$ the $\eta$-locus of $f$ in $M$ is the transverse preimage $\hat\kappa_f^{-1}(B_\eta(n))$. This observation was originally proved by Thom (see \cite{thomthm}); see also \cite{sasohm} for a proof specialised to finitely determined singularities. In the proof we use the maps that induce the tangent bundles $TM$ and $TP$ from $\gamma_n$ and $\gamma_{n+k}$ respectively, hence this cannot be directly applied to quasi-holomorphic maps between real smooth manifolds. However, 
it will turn out that a ``stabilised'' version of the Kazarian space has a property for quasi-holomorphic maps analogous to the above one.

\begin{lemma}
  \label{lemma:cplxnb}
  If $f\colon M^n\to P^{n+2k}$ is a quasi-holomorphic map, then the stable normal bundle $\nu_f:=f^*TP\ominus TM$ is endowed with a complex structure.
\end{lemma}

\begin{prf}
  A vector bundle representing $\nu_f$ is $f^*TP\oplus\nu_M$, where $\nu_M$ is defined as the normal bundle of an embedding $e\colon M\into\C^L$ where $L$ is such a large number that the embeddings $M\into\C^L$ are all isotopic. Then $f^*TP\oplus\nu_M$ is the normal bundle of the embedding $f\times e\colon M\into P\times\C^L$.

  Recall from definition \ref{defi:qhol} that $f$ maps the strata $\sigma_i$ in $M$ to (the closures of the) strata $\tilde\sigma_{i+k}$ in $P$ and let $\nu_i$ and $\tilde\nu_{i+k}$ denote the normal bundles of these strata in $M$ and $P$ respectively (with $\tilde\nu_{i+k}$ naturally extended to a bundle along the immersion $f|_{\sigma_i}$). Then $\nu_i$ is a subbundle of the pullback $f|_{\sigma_i}^*\tilde\nu_{i+k}\oplus\varepsilon^L$ (where $\varepsilon^L$ is the trivial complex $L$-vector bundle), and the normal bundle of the restriction of the embedding $f\times e$ to $\sigma_i$ is $\nu_i\oplus(f|_{\sigma_i}^*\tilde\nu_{i+k}\oplus\varepsilon^L)/\nu_i$. Thus the normal bundle of the whole embedding $f\times e$ is composed of the vector bundles $(f|_{\sigma_i}^*\tilde\nu_{i+k}\oplus\varepsilon^L)/\nu_i$ over the strata $\sigma_i$, which are all endowed with complex structures by definition. Moreover, the complex structures given on different strata fit together to a complex structure on the whole normal bundle because of the compatibility condition on the structures $J_i$ and $\tilde J_{i+k}$ for different $i$'s in definition \ref{defi:cstrat}. Thus the representatives of the stable vector bundle $\nu_f$ are endowed with complex structures in a natural way, and so $\nu_f$ is also complex.
\end{prf}

The space we introduced in definition \ref{defi:kaz} is the ``unstable'' version of the Kazarian space; its stabilisation will involve the natural suspension inclusions $J_0^\infty(\C^n,\C^{n+k})\into J_0^\infty(\C^{n+1},\C^{n+k+1})$. These suspensions do not behave nicely with respect to singularities that are finitely determined but not stable (see \cite{sasohm}), hence in the rest of this section $\sigma$ will be a set of stable singularities. Now instead of the jet bundle $J_0^\infty(\gamma_n,\gamma_{n+k})\to B\U(n)\times B\U(n+k)$ we will consider $J_0^\infty(\varepsilon^n,\gamma_{n+k})\to B\U(n+k)$ (where, again, $\varepsilon^n$ is the trivial complex $n$-vector bundle). 
We will denote by $K_\sigma(n)\subset J_0^\infty(\varepsilon^n,\gamma_{n+k})$ the subspace we get by taking in each fibre the subspace $J_\sigma(n)\subset J_0^\infty(\C^n,\C^{n+k})$ and note that since the above suspension inclusion maps $J_\sigma(n)$ into $J_\sigma(n+1)$, its natural extension to jet bundles maps $K_\sigma(n)$ into $K_\sigma(n+1)$.

\begin{defi}
  \label{defi:stkaz}
  For a set $\sigma$ of singularities of stable holomorphic germs, the \textit{stable Kazarian space} of $\sigma$-maps is the direct limit $K_\sigma:=\underset{n\to\infty}\lim K_\sigma(n)$.
\end{defi}

Similarly to the unstable Kazarian space, we also have that $K_\sigma$ is the union of subspaces homotopy equivalent to $BG_\eta$ (for all $[\eta]\in\sigma$) and the normal bundle of $BG_\eta$ in $K_\sigma$ is $\xi_\eta$. Now the following theorem can be proved in completely the same way as its unstable analogue for holomorphic maps (see e.g. \cite{thomthm}, \cite{kaztp}), hence we will only sketch the proof.

\begin{thm}
  \label{thm:kaz}
  For a set $\sigma$ of singularities of stable holomorphic germs, if $f\colon M^n\to P^{n+2k}$ is a quasi-holomorphic $\sigma$-map, then there is a map $\kappa_f\colon M\to K_\sigma$ such that for all singularities $[\eta]\in\sigma$ we have $\kappa_f^{-1}(BG_\eta)=\eta(f)$. For almost all maps $f$ the map $\kappa_f$ is transverse to the closure of $BG_\eta$.
\end{thm}

\begin{sketchprf}
  Fix Riemannian metrics on $M$ and $P$ which are compatible with the complex structures $J_i,\tilde J_i$ on the normal bundles of the strata $\sigma_i,\tilde\sigma_i$ given by the quasi-holomorphic structure (for all $i$); see definition \ref{defi:qhol}. Denote by $\exp_M$ and $\exp_P$ the exponential maps on $M$ and $P$ respectively and let $\overset{_\circ}{T}M\subset TM$ be a neighbourhood of the zero-section such that $\exp_M|_{{\overset{_\circ}{T}}_pM}$ is a diffeomorphism onto its image for all $p\in M$. Then there is a fibrewise map $F$ that makes the following diagram commutative:
  $$\xymatrix{
    \overset{_\circ}{T}M\ar@{-->}[r]^{F}\ar[d]_{\exp_M} & TP\ar[d]^{\exp_P} \\
    M\ar[r]^f & P.
  }$$
  Let $L$ be such a large number that $M$ can be embedded into $\C^L$ uniquely up to isotopy and denote by $\nu_M$ the (uniquely determined) normal bundle of $M\into\C^L$. If $J_0^\infty(TM\oplus\nu_M,f^*TP\oplus\nu_M)=J_0^\infty(\varepsilon^L,f^*TP\oplus\nu_M)\to M$ is the bundle of infinite jets of germs of fibrewise maps $\varepsilon^L\to f^*TP\oplus\nu_M$ along the zero-section, then the map $F\oplus\id_{\nu_M}$ defines a section $\tilde F\colon M\to J_0^\infty(\varepsilon^L,f^*TP\oplus\nu_M)$ of this jet bundle. Note that the fibre of the bundle $J_0^\infty(\varepsilon^L,f^*TP\oplus\nu_M)$ is $J_0^\infty(\C^L,\C^{L+k})$ 
  and $\tilde F$ maps into the subbundle $J_\sigma(M,L)\subset J_0^\infty(\varepsilon^L,f^*TP\oplus\nu_M)$ formed by the fibres $J_\sigma(L)$ since $f$ is a $\sigma$-map.

  The fibre bundle $J_\sigma(M,L)\xra{J_\sigma(L)}M$ is induced from the universal bundle with fibre $J_\sigma(L)$, that is, from the restriction $K_\sigma(L)\xra{J_\sigma(L)}B\U(L+k)$ of the jet bundle $J_0^\infty(\varepsilon^L,\gamma_{L+k})\to B\U(L+k)$, by a map $M\to B\U(L+k)$. This map also induces the vector bundle $f^*TP\oplus\nu_M$ from $\gamma_{L+k}$, thus its composition with the inclusion $B\U(L+k)\subset B\U$ is the map $M\to B\U$ inducing the stable normal bundle $\nu_f$. Hence we obtain a pullback square
  $$\xymatrix{
    J_\sigma(M,L)\ar[r]^{\kappa}\ar[d]_{J_\sigma(L)} & K_\sigma(L)\ar[d]^{J_\sigma(L)} \\
    M\ar@/^4pc/[u]^{\tilde F}\ar[r] & B\U(L+k)
  }$$
  with the section $\tilde F$ defined above. Then the map $\kappa_f$ is defined to be the composition of $\kappa\circ\tilde F$ with the inclusion $K_\sigma(L)\subset K_\sigma$. The preimage of the subset $BG_\eta\subset K_\sigma$ under the map $\kappa_f$ is now the preimage of the corresponding $\eta$-locus under the map $\tilde F$ which is the $\exp_M$-preimage of the $\eta$-locus of $F$, that is, $\eta(f)$. The claimed transversality of $\kappa_f$ to $\ol{BG}_\eta$ follows from usual transversality properties. 
\end{sketchprf}

\begin{rmk}
  In particular, the above theorem shows that the global normal form of a monosingularity $[\eta]$ also has the property that the normal bundle of the $\eta$-locus in the source manifold of any quasi-holomorphic map (in the complement of the higher codimensional strata) is induced from the bundle $\xi_\eta\to BG_\eta$; cf. lemma \ref{lemma:semig}.
\end{rmk}

\section{Global singularity theory of quasi-holomorphic maps}
\label{sec:globsingth}

The universal map $f_\tau\colon Y_\tau\to X_\tau$ and the Kazarian space $K_\sigma$ described in the previous section are natural tools for the study of global singularity properties of quasi-holomorphic maps, as theorems \ref{thm:univmap} and \ref{thm:kaz} show. In the following two subsections we specialise two central elements of global singularity theory to quasi-holomor\-phic maps, and show basic properties and characteristics of them.

\subsection{(Multi)singular loci and Thom polynomials}
\label{ssec:loci}

For a multisingularity set $\tau$, the construction of the map $f_\tau\colon Y_\tau\to X_\tau$ immediately implies that for each multisingularity $\mu\in\tau$ and any $\tau$-map $f\colon M\to P$ the locus $\mu(f)$ is the preimage of the subspace $B_\mu\subset Y_\tau$. This means that the properties of the space $Y_\tau$ strongly correspond to the multisingular loci of $\tau$-maps. These properties were analysed and applied by Rimányi \cite{rim} to the computation of Thom polynomials, incidences of singularities and multiple point formulae for holomorphic maps. Similarly, if $\sigma$ is a singularity set without multiplicities, then the properties of the space $K_\sigma$ correspond to the singular loci of $\sigma$-maps. This again has been applied by many to the investigation of Thom polynomials, avoiding ideals and related objects for holomorphic maps; see e.g. \cite{kaztp} for a review on these methods. Now theorems \ref{thm:univmap} and \ref{thm:kaz} imply that all properties of the (multi)singular loci of holomorphic maps deduced from investigating these spaces are the exact same for those of quasi-holomorphic maps. Here we will only highlight this statement on Thom polynomials, which can be considered as the most important invariants of singular loci.

Recall that for any stable singularity $[\eta]$ of $k$-codimensional holomorphic germs, there is a fixed element $\Tp_\eta\in\Z[c_1,c_2,\ldots]$, called the \textit{Thom polynomial} of $[\eta]$, such that for any stable holomorphic map $f\colon M^n\to P^{n+k}$ the cohomology class in $H^*(M;\Z)$ represented by the closure of the $\eta$-locus of $f$ (i.e. the Poincaré dual of $\ol\eta(f)$) is the evaluation $\Tp_\eta(f)$ defined by substituting the normal Chern class $c_i(\nu_f)$ in the variable $c_i$ of $\Tp_\eta$ for each $i$; see e.g. \cite{tp}, \cite{kaztp}. The Thom polynomial $\Tp_\eta$ is usually defined as an element in $H^*(B\U;\Z)\cong\Z[c_1,c_2,\ldots]$ and the evaluation $\Tp_\eta(f)$ is the pullback of $\Tp_\eta$ by the inducing map of the stable normal bundle $\nu_f$.

A more concrete description of the cohomology class $\Tp_\eta\in H^*(B\U;\Z)$ is the following. Let $\sigma$ be the set of all stable singularities of $k$-codimensional holomorphic germs, and denote by $p\colon K_\sigma\to B\U$ the map induced by the restrictions $K_\sigma(n)\to B\U(n+k)$ of the bundle maps $J_0^\infty(\varepsilon^n,\gamma_{n+k})\to B\U(n+k)$. The subvariety $\ol{BG}_\eta\subset K_\sigma$ determines a Poincaré dual cohomology class in $H^*(K_\sigma;\Z)$. The Thom polynomial $\Tp_\eta$ is the image of this class under the pullback $p^*$. Now theorem \ref{thm:kaz} immediately implies the following:

\begin{thm}
  \label{thm:tp}
  Let $[\eta]$ be a stable holomorphic singularity and $f\colon M\to P$ a stable quasi-holomorphic map. Then the Poincaré dual cohomology class of the closure $\ol\eta(f)$ of the $\eta$-locus in $H^*(M;\Z)$ exists, and this class is $\Tp_\eta(f)$.
\end{thm}

%

\begin{ex}
  The Thom polynomial of the fold singularity of $k$-codimensional maps is $\Tp_{A_1}=c_{k+1}$. As we remarked before, a double branched cover is a quasi-holomorphic fold map of codimension $0$. Thus for a double branched cover $f\colon M^n\to P^n$, the cohomology class in $H^*(M;\Z)$ represented by the branch locus is $c_1(\nu_f)=f^*c_1(P)-c_1(M)$.
\end{ex}

\begin{rmk}
  One can also define the \textit{restricted Thom polynomial} of a singularity $[\eta]$ with respect to any fixed singularity set $\sigma$ containing $[\eta]$; see e.g. \cite{rim}. This can be given by taking the dual cohomology class of $\ol{BG}_\eta$ in $K_\sigma$, and similarly to theorem \ref{thm:tp} we get that for any $\sigma$-map $f\colon M\to P$, its pullback by the map $\kappa_f$ is the cohomology class represented by $\ol\eta(f)$ in $H(M;\Z)$. We note that for example the Thom polynomials of Morin singularities $A_r$ are not known in general, but the restricted Thom polynomial of $A_r$ to the singularity set $\Sigma^0\cup\{A_i\mid i=1,2,\ldots\}$ (i.e. to Morin maps) is known for all $r$; see \cite{morchar}.
\end{rmk}

\begin{rmk}
  Thom polynomials exist for not only stable but also finitely determined singularities in general; see \cite{sasohm} where some of these polynomials were computed. However, the Thom polynomials of unstable singularities cannot be written in Chern classes of the stable normal bundle, instead they are polynomials in the tangent Chern classes of the source and target manifolds. This makes sense for holomorphic maps but cannot be extended to quasi-holomorphic maps, reflecting the fact that the stable Kazarian space is only defined for sets of stable singularities.
\end{rmk}

\subsection{Cobordism groups}
\label{ssec:cob}

After considering $Y_\tau$ and $K_\sigma$, which ``classified'' the (multi)singular loci of $\tau$-maps and $\sigma$-maps respectively, let us now turn to the objects classified by the space $X_\tau$. For the rest of this section we fix the set $\tau$ of multisingularities of finitely determined holomorphic germs of codimension $k$.

\begin{defi}
  \label{defi:cob}
  Let $P^{n+2k}$ be a manifold without boundary. Two $\tau$-maps $f_0\colon M_0^n\to P^{n+2k}$ and $f_1\colon M_1^n\to P^{n+2k}$ with closed source manifolds $M_0$ and $M_1$ are said to be \textit{cobordant} if there is
  \begin{enumerate}[label=(\roman*)]
  \item a compact manifold $W^{n+1}$ with boundary $\partial W=M_0\sqcup M_1$,
  \item a $\tau$-map $F\colon W\to P\times[0,1]$ transverse to the boundary such that for $i=0,1$ we have $F^{-1}(P\times\{i\})=M_i$ and the restriction $F|_{M_i}\colon M_i\to P\times\{i\}$ is the map $f_i$.
  \end{enumerate}
  Cobordism is an equivalence relation on the set of $\tau$-maps of $n$-manifolds to $P$. The cobordism class of $f\colon M\to P$ is denoted by $[f]$ and the set of these cobordism classes is denoted by $\Qhol_\tau(P)$; in the case $P=\R^{n+2k}$ we simplify this notation to $\Qhol_\tau(n,k)$.
\end{defi}

If the multisingularity set $\tau$ is complete, then $\Qhol_\tau(P)$ admits a natural commutative semigroup operation: we define the operation applied to the cobordism classes of two $\tau$-maps $f\colon M\to P$ and $g\colon N\to P$ by $[f]+[g]:=[f\sqcup g]$ where the disjoint union is $f\sqcup g\colon M\sqcup N\to P$. The neutral element of this operation is represented by the empty map. If the manifold $P$ is of the form $Q\times\R^1$, then this is an Abelian group operation where the inverse of $[f]$ is defined by composing $f$ with the reflection $\rho$ with respect to $Q\times\{0\}$ in the space $Q\times\R$; this is justified by the null-cobordism of $f\sqcup\rho\circ f$ given by the rotation in $Q\times\R\times\R_+$ around the axis $Q\times\{(0,0)\}$. Moreover, if we have again $P=Q\times\R$ but $\tau$ is not complete, then this group operation still exists if we define $f\sqcup g$ as the ``far away'' disjoint union given by composing $g$ with an appropriate translation in $Q\times\R$ along $\R$. However, in this case the operation is not necessarily Abelian unless $P$ is of the form $Q'\times\R^2$. Later on, when talking about the set $\Qhol_\tau(P)$, we will tacitly endow it with the above described algebraic structures, hence in most cases we will call it a \textit{cobordism group}.

\begin{ex}
  Note that several classical examples of cobordism groups of maps are defined analogously to definition \ref{defi:cob}. In particular we have the following:
  \begin{enumerate}
  \item For any vector bundle $\zeta$ and manifold $P$ there is a cobordism set (or, in many cases, group) $\Emb^\zeta(P)$ of embeddings $M\into P$ with normal bundle induced from $\zeta$. The usual Pontryagin--Thom construction shows that $\Emb^\zeta(P)$ is isomorphic to $[\cpt P,T\zeta]$ 
    (in particular, if $P$ is a Eucledian space, then this is a homotopy group of $T\zeta$). In other words, $T\zeta$ is the \textit{classifying space} of cobordisms of such embeddings.
  \item Similarly, if $\Imm^\zeta(P)$ denotes the cobordism group of immersions $M\imto P$ with normal bundle induced from $\zeta$, then we have $\Imm^\zeta(P)\cong[S^\infty\cpt P,S^\infty T\zeta]\cong[\cpt P,\Gamma T\zeta]$ (recall that $\Gamma=\Omega^\infty S^\infty$). 
    That is, $\Gamma T\zeta$ is the \textit{classifying space} of cobordisms of immersions (see \cite{wells}).
  \item Again similarly one can define the cobordism set (or group) $\Imm_r^\zeta(P)$ of immersions as above but without $(r+1)$-tuple points; their classifying space was constructed in \cite{limm} and will be denoted by $\Gamma_r(\zeta)$. 
    We note that earlier Uchida introduced another, related notion of cobordisms of immersions with restricted multiplicities in \cite{eximm} but did not construct their classifying space.
  \item Classifying spaces of cobordism groups of real smooth maps with various types of (multi)singu\-larity restrictions were constructed by many; most notably see e.g. \cite{rsz}, \cite{hosszu} and the different constructions of Ando \cite{ando} and Sadykov \cite{sad}. For us the main example is the real smooth version of the space $X_\tau$ in \cite{rsz}.
  \item A similar construction for the classifying space of cobordisms of simple branched coverings was given in \cite{nagy}.
  \end{enumerate}
\end{ex}


Now in the case of quasi-holomorphic maps with restricted multisingularities we have the following classifying property of the space $X_\tau$. Its proof is the complete analogue of that of the main theorem in \cite{rsz}, and so we omit it here. 

\begin{thm}
  \label{thm:cob}
  For any set $\tau$ of multisingularities of finitely determined holomorphic germs and any manifold $P$ we have
  $$\Qhol_\tau(P)\cong[\cpt P,X_\tau]$$
\end{thm}

%

\begin{rmk}
  The classifying space of cobordisms of $\tau$-maps is unique up to homotopy because the functor $\Qhol_\tau$ is Brown representable, and $X_\tau$ is a specifically constructed representative of this homotopy type. Of course, this statement is not precise, since to apply the Brown representability theorem we first have to extend $\Qhol_\tau$ to a contravariant functor of simplicial complexes. We will not do this here, but refer to \cite[part I]{hosszu}, where this was clarified for cobordisms of real smooth maps with prescribed singularities; all methods there also work in the quasi-holomorphic setting.
\end{rmk}

\begin{exenum}
\item If $\tau$ is the multisingularity set $\{i\Sigma^0\mid i=1,\ldots,r\}$ for some $1\le r\le\infty$, then $X_\tau$ is the classifying space of cobordisms of immersions with at most $r$-tuple points and complex normal bundles. In particular, for $r=1$ we have $X_\tau\cong T\gamma_k$ and for $r=\infty$ we have $X_\tau\cong\Gamma T\gamma_k$.
\item If the codimension $k$ is $0$ and $\tau$ is any decreasing subset of the multisingularity set $\N\la r\Sigma^0,$ $1[(-)^r]\ra$, then the space $X_\tau$ is the classifying space of $r$-sheeted branched coverings with a corresponding multiplicity restriction. In particular, for $\tau=\{r\Sigma^0,1[(-)^r]\}$ the space $X_\tau$ is homotopy equivalent to the one constructed in \cite{nagy}.
\end{exenum}

\section{Computations on cobordism groups}
\label{sec:cob}

As theorem \ref{thm:cob} shows, the computation of cobordism groups of quasi-holomorphic maps corresponds to the computation of the groups of homotopy classes of maps to the classifying space; in particular we have $\Qhol_\tau(n,k)\cong\pi_{n+2k}(X_\tau)$. In the following we will construct and analyse exact sequences that contain such homotopy groups. For simplicity of notation, if $\tau=\N\la\sigma\ra$ is a complete multisingularity set (for a singularity set $\sigma$), then we will put $X_\sigma:=X_{\N\la\sigma\ra}$ and $\Qhol_\sigma:=\Qhol_{\N\la\sigma\ra}$.

\subsection{The key fibration}
\label{ssec:key}

In \cite{hosszu} a fibration of classifying spaces of cobordism groups of (real smooth) singular maps was constructed if the codimension of the maps is positive, yielding exact sequences of cobordism groups as homotopy long exact sequences. This fibration was later directly obtained by Terpai \cite{key} in a somewhat more general form. The proof in \cite{key} goes through in a completely analogous way if we consider quasi-holomorphic maps instead of real smooth ones, showing the following theorem. We will again only sketch the proof. 

\begin{thm}
  \label{thm:key}
  Let $\tau$ be a set of multisingularities of finitely determined holomorphic germs of complex codimension $k>0$, and suppose that $\tau$ is of the form $\{i[\eta]+\mu\mid i=1,\ldots,r;\mu\in\tau'\}$ where $\tau'$ is a complete multisingularity set which contains all singularities below $[\eta]$ (and does not contain $[\eta]$) and $1\le r\le\infty$. Then there is a Serre fibration
  $$X_\tau\xra{X_{\tau'}}\Gamma_r(\tilde\xi_\eta).$$
\end{thm}

\begin{sketchprf}
  For a $\tau$-map $f\colon M^n\to P^{n+2k}$, the restriction $f|_{\eta(f)}$ is an immersion $\eta(f)\imto P$ without $(r+1)$-tuple points and with normal bundle induced from $\tilde\xi_\eta$. Sending the cobordism class of $f$ as a $\tau$-map to the cobordism class of $f|_{\eta(f)}$ as an immersion with this structure defines a natural transformation of functors $\Qhol_\tau\Rightarrow\Imm^{\tilde\xi_\eta}_r$. A version of the Brown representability theorem now implies that this transfromation is represented by a map $p\colon X_\tau\to\Gamma_r(\tilde\xi_\eta)$ between the corresponding classifying spaces (actually this version is an immediate consequence of the representability theorem and the Yoneda lemma). We note here that the map $p$ could also be obtained directly by referring to the constructions of the spaces $X_\tau$ in theorem \ref{thm:univmap} and $\Gamma_r(\tilde\xi_\eta)$ in \cite{limm}. The rest of the proof consists of two steps: firstly, proving that $p$ is a Serre fibration, and secondly, proving that the fibre of $p$ is homotopy equivalent to $X_{\tau'}$.

  The first step is proved by directly checking the homotopy lifting property of the map $p$. This, when translated to the cobordisms classified by the spaces $X_\tau$ and $\Gamma_r(\tilde\xi_\eta)$, means that we are given a manifold $P^{n+2k}$, a $\tau$-map $f\colon M^n\to P=P\times\{0\}$ and an immersion $i\colon N^{n-c_\eta+1}\imto P\times[0,1]$, which maps the boundary $\partial N$ to the boundary $P\times\{0,1\}$ trasversely, has no $(r+1)$-tuple points and its normal bundle is induced from $\tilde\xi_\eta$, moreover, we also have $i^{-1}(P\times\{0\})=\eta(f)$ and $i|_{\eta(f)}=f|_{\eta(f)}$. Our goal is to find a $\tau$-map $F\colon W^{n+1}\to P\times[0,1]$ that extends $f\colon M\to P\times\{0\}$, maps $\partial W\setminus M$ to $P\times\{1\}$, is transverse to the boundary, and for which $\eta(F)=N$ and $F|_{\eta(f)}=i$.


  Now the normal structure of the immersion $i$ and lemma \ref{lemma:pb} allow us to extend $i$ to a quasi-holomorphic map $F_\eta$ of a disk bundle $U$ over $N$ (which is induced from the 
  disk bundle 
  $D\xi_\eta$). 
  The boundary of $U$ is the union of a sphere bundle over $N$, denoted by $\partial_SU$, and two disk bundles over $i^{-1}(P\times\{0\})$ and $i^{-1}(P\times\{1\})$, denoted by $\partial_0U$ and $\partial_1U$ respectively. By lemma \ref{lemma:semig} we can assume that $\partial_0U$ is a closed tubular neighbourhood of $\eta(f)$ and the restrictions $F_\eta|_{\partial_0U}$ and $f|_{\partial_0U}$ coincide. We also extend the map $f\colon M\to P\times\{0\}$ by a collar to the map $f\times\id_{[-1,0]}\colon M\times[-1,0]\to P\times[-1,0]$. Then we can glue $U$ to $M\times[-1,0]$ at $\partial_0U\subset M\times\{0\}$ resulting (after smoothing corners) in an $(n+1)$-manifold $V_0$ with boundary the union of $M_0:=M\times\{-1\}$ and a manifold $M_1$ composed of $M\times\{0\}\setminus\partial_0U$ and $\partial_SU\cup\partial_1U$. The union of the maps $f\times\id_{[-1,0]}$ and $F_\eta$ is now a mapping $F_0$ of $V_0$ to $P\times[-1,1]\approx P\times[0,1]$ that satisfies almost every property that we expect from the proposed $\tau$-map $F\colon W\to P\times[0,1]$; the only exception is that it does not map the boundary part $M_1$ to the boundary part $P\times\{1\}$. To solve this problem, we choose a nowhere vanishing ``outward'' normal vector field $v$ at $\partial_SU\subset M_1$. By the ``stratified compression theorem'' \cite[theorem 1]{hosszu} (cf. \cite{rs}) we may assume (after a small ambient isotopy in $P\times[0,1]$) that $v$ is vertical, i.e. everywhere in the subbundle $T[0,1]<T(P\times[0,1])$. Then we can define another $(n+1)$-manifold
  $$V_1:=\{(q,t)\in(M_1\setminus\partial_1U)\times[0,1]\mid F_0(q)=(q',s)\in P\times[0,1],t\in[s,1]\}$$
  and put $F_1\colon V_1\to P\times[0,1]$ to be the mapping $(q,t)\mapsto(F_0(q),t)$. The boundary of $V_1$ is the union of $M_1\setminus\partial_1U$ and a manifold $M_2$ mapped to $P\times\{1\}$. Then the manifold $W^{n+1}$ formed by $V_0$ and $V_1$ glued together at $M_1\setminus\partial_1U$ and the map $F\colon W\to P\times[0,1]$ defined as the union $F_0\cup F_1$ is such that its smoothing satisfies every property we required. This concludes the first step.
  
  The second step is almost identical to the first one. The fibre of the map $p\colon X_\tau\to\Gamma_r(\tilde\xi_\eta)$ is originally the classifying space of cobordisms of $\tau$-maps whose restriction to the $\eta$-locus is null-cobordant; we want to prove that it is also the classifying space of cobordisms of $\tau'$-maps. Again we translate this problem to the cobordisms classified by the spaces $X_\tau$, $X_{\tau'}$ and $\Gamma_r(\tilde\xi_\eta)$. Now we are given a $\tau$-map $f\colon M\to P\times\{0\}$ and an immersion $i\colon N\imto P\times[0,1]$ as above, the only difference now is that $i$ maps the whole of the boundary $\partial N$ to $P\times\{0\}$, i.e. $i$ is a null-cobordism of $f|_{\eta(f)}$. Again we want to extend $f$ to a $\tau$-map $F\colon W^{n+1}\to P\times[0,1]$ as above, with $\eta(f)=N$ and $F|_{\eta(F)}=i$, that is, $F$ should be a cobordism between $f$ and a $\tau'$-map. The construction of $F$ is exactly the same as that in the first step. This construction indeed proves that the fibre of the map $p$ is the classifying space of cobordisms of $\tau'$-maps, hence by the Brown representability theorem this fibre is homotopy equivalent to $X_{\tau'}$.
\end{sketchprf}

\begin{rmk}
  \label{rmk:key}
  In the proof we used perturbations of maps that only change them at specific parts, keeping the rest fixed (see the smoothings and the application of the compression theorem). These do not exist for holomorphic maps, but they work for quasi-holomorphic ones. We also note that when applying the compression theorem, we 
  used the codimension assumption $k>0$ since the theorem only holds for maps of codimension at least $2$ and we applied it to the mapping of an $n$-manifold into an ($n+2k+1$)-manifold.
\end{rmk}


The above theorem is a key tool for computing cobordism groups of $\tau$-maps and it is especially useful in the case of complete sets of multisingularities $\tau$ (i.e. when $r=\infty$), that is, when there are no global restrictions, only local ones. In this case the base space $\Gamma_\infty(\tilde\xi_\eta)$ is the classifying space of immersions with normal bundle induced from $\tilde\xi_\eta$ which is $\Gamma T\tilde\xi_\eta=\Omega^\infty S^\infty T\tilde\xi_\eta$. Now if $\sigma$ is a set of holomorphic singularities and $\eta\in\sigma$ is a maximal element, then we will write
$$p_\sigma^\eta\colon X_\sigma\xra{X_{\sigma\setminus\{\eta\}}}\Gamma T\tilde\xi_\eta$$
for the fibration in the theorem and will call it the \textit{key fibration}.

\subsection{Rational cobordism groups}
\label{ssec:rat}

In the following we will show that the rationalisations of the cobordism groups of 
$\sigma$-maps can be easily understood for any set $\sigma$ 
of 
singularities without multiplicities. 
We will denote by $[-]_\Q$ the rational homotopy type functor. Before stating our theorem let us recall the following straightforward corollary of \cite[lemma 81]{hosszu} and the Dold--Thom theorem:

\begin{lemma}
  \label{lemma:dt}
  For any connected space $X$ there is a homotopy equivalence
  $$[\Gamma X]_\Q\cong\prod_{i=1}^\infty K(H_i(X;\Q),i),$$
  in particular, the stable Hurewicz homomorphism is a rational isomorphism
  $$\pi^s_*(X)\otimes\Q\cong\pi_*(\Gamma X)\otimes\Q\cong H_*(X;\Q).$$
\end{lemma}

Now we can give a rational description of the classifying space $X_\sigma$. We note that this description has been stated in \cite[claim 121]{hosszu} as a special case of a more general observation; the proof we give here is essentially the same as the one given there. We also note however, that \cite{hosszu} lacks the precise introduction to quasi-holomorphic maps and the classifying spaces of their cobordism groups, hence the proof there is incomplete.

\begin{thm}
  \label{thm:keyractrivi}
  For any set $\sigma$ of singularities of finitely determined positive codimensional holomorphic germs we have 
  $$[X_\sigma]_\Q\cong\prod_{\eta\in\sigma}[\Gamma T\tilde\xi_\eta]_\Q.$$
\end{thm}

\begin{prf}
  We proceed by induction on an arbitrary complete order extending the natural partial order of the singularities in $\sigma$. The base case is obvious since the lowest element is the ``non-singular singularity'' $\Sigma^0$ and $X_{\{\Sigma^0\}}$ is just $\Gamma T\tilde\xi_{\Sigma^0}=\Gamma T\gamma_k$. For the induction step let $\eta\in\sigma$ be the maximal element and put $\sigma'=\sigma\setminus\{\eta\}$; suppose that the theorem is satisfied for $\sigma'$, i.e. we have
  $$[X_{\sigma'}]_\Q\cong\prod_{\vartheta\in\sigma'}[\Gamma T\tilde\xi_\vartheta]_\Q.$$

  For any singularity $\vartheta$ the vector bundle $\tilde\xi_\vartheta$ is a complex bundle (hence of even real rank) over $BG_\vartheta$ where $G_\vartheta$ is a compact $\C$-linear group. Hence the odd homology groups of the space $BG_\vartheta$ vanish, and so, by the Thom isomorphism, so do those of $T\tilde\xi_\vartheta$. Thus, using lemma \ref{lemma:dt}, we have
  $$[\Gamma T\tilde\xi_\vartheta]_\Q\cong\prod_{i=1}^\infty K(H_{2i}(T\tilde\xi_\vartheta;\Q),2i).$$
  Applying this 
  for all $\vartheta\in\sigma'$ 
  yields that 
  $[X_{\sigma'}]_\Q$ is a product of rational Eilenberg--MacLane spaces of even indices.

  The rationalisation of the key fibration $p_\sigma^\eta$ is of the form
  $$[X_\sigma]_\Q\xra{[X_{\sigma'}]_\Q}[\Gamma T\tilde\xi_\eta]_\Q$$
  (see e.g. \cite[5.1]{bk}). This is induced from the universal fibration with fibre $[X_{\sigma'}]_\Q$ by a map $b$, i.e. there is a pullback diagram
  $$\xymatrix{
    [X_\sigma]_\Q\ar[r]\ar[d]_{[X_{\sigma'}]_\Q} & \mathrm{*}\ar[d]^{[X_{\sigma'}]_\Q} \\
    [\Gamma T\tilde\xi_\eta]_\Q\ar[r]^{b} & B[X_{\sigma'}]_\Q.
  }$$
  Here the base space $B[X_{\sigma'}]_\Q$ is a product of Eilenberg--MacLane spaces of odd indices:
  $$B\left(\prod_{\vartheta\in\sigma'}\prod_{i=1}^\infty K(H_{2i}(T\tilde\xi_\vartheta;\Q),2i)\right)\cong\prod_{\vartheta\in\sigma'}\prod_{i=1}^\infty BK(H_{2i}(T\tilde\xi_\vartheta;\Q),2i)\cong\prod_{\vartheta\in\sigma'}\prod_{i=1}^\infty K(H_{2i}(T\tilde\xi_\vartheta;\Q),2i+1),$$
  while $[\Gamma T\tilde\xi_\eta]_\Q$ is a product of ones of even indices. Thus any map $[\Gamma T\tilde\xi_\eta]_\Q\to B[X_{\sigma'}]_\Q$ induces the zero homomorphism in cohomology, in particular $b^*$ maps each fundamental class to zero, hence $b$ is null-homotopic. This means that the key fibration $p_\sigma^\eta$ is rationally trivial, i.e. we have
  $$[X_{\sigma}]_\Q\cong[X_{\sigma'}]_\Q\times[\Gamma T\tilde\xi_\eta]_\Q.$$
  This concludes the induction step.
\end{prf}

\begin{rmk}
  Geometrically the above theorem means that up to rational cobordism, a quasi-holomorphic map behaves the same as the collection of its restrictions to the individual singular loci as immersions independent of each other.
\end{rmk}

From theorem \ref{thm:keyractrivi} we immediately obtain the following corollaries, in all of which the singularity set $\sigma$ is as in the theorem.

\begin{crly}
  \label{crly:ractrivi}
  For any 
  manifold $P^{n+2k}$ we have
  $$\Qhol_\sigma(P)\otimes\Q\cong\bigoplus_{\eta\in\sigma}\Imm^{\tilde\xi_\eta}(P)\otimes\Q.$$
\end{crly}

In the case $P^{n+2k}=\R^{n+2k}$ this specialises to the following (using lemma \ref{lemma:dt} and the Thom isomorphism):

\begin{crly}
  \label{crly:rachol}
  We have
  $$\Qhol_\sigma(n,k)\otimes\Q\cong\bigoplus_{\eta\in\sigma}H_{n-2c_\eta}(BG_\eta;\Q).$$
\end{crly}

Since $G_\eta$ is a compact $\C$-linear group for any singularity $\eta$, the classifying space $BG_\eta$ has no odd homology. Thus we also get the following:

\begin{crly}
  If $n$ is odd, then the group $\Qhol_\sigma(n,k)$ is torsion.
\end{crly}


\begin{ex}
  Denote by 
  $\mu_r$ the set of the Morin singularities $A_i$ for $i\le r$ (we also allow $r=\infty$). For any $i$ the codimension of the singularity $A_i$ is $i(k+1)$ and the group $G_{A_i}$ is $\U(1)\times\U(k)$ for $i\ge1$ (see e.g. \cite{rim}) and $\U(k)$ for $i=0$. Hence by corollary \ref{crly:rachol} we have
  \begin{alignat*}2
    \Qhol_{\mu_r}(n,k)\otimes\Q&\cong H_n(B\U(k);\Q)\oplus\bigoplus_{i=1}^rH_{n-2i(k+1)}(B\U(1)\times B\U(k);\Q)\cong\\
    &\cong\Q^{p_k\left(\frac{n}2\right)+\sum_{i=1}^r\sum_{j=0}^{\frac{n}2-i(k+1)}p_k\left(\frac{n}2-j\right)}
  \end{alignat*}
  where $p_k(m)$ denotes the number of partitions of an integer $m$ to sums of positive integers not greater than $k$.
\end{ex}

\subsection{2-codimensional fold cobordism groups}
\label{ssec:fold}

In the previous section we saw that 
the free part (i.e. the factor by the torsion part) of the cobordism groups of quasi-holomorphic maps can be completely determined in all cases (up to the computation of maximal compact subgroups of symmetry groups of singularities). The determination of the torsion part is considerably harder, as we shall now demonstrate by obtaining further results on these 
cobordism groups in the case $k=1$ and $\sigma=\mu_1=\{\Sigma^0,A_1\}$, i.e. for fold maps of (real) codimension $2$, which can be considered as the simplest (non-regular) case. 

\begin{rmk}
  If we also forbid the fold singularity, that is, we only consider 2-codimen\-sional immersions with complex structures on their normal bundles, then the classifying space we obtain is $\Gamma T\gamma_1\cong\Gamma\CP^\infty$, hence the cobordism group of quasi-holomorphic immersions of $n$-manifolds into $\R^{n+2}$ is $\pi^s_{n+2}(\CP^\infty)$. This group was determined for indices $n\le16$ by Liulevicius \cite{liu}, Mosher \cite{mosh} and Mukai \cite{muk} by quite involved algebraic computations. This suggests that exact computations in the case of fold maps are even less feasible, however, we will see that at least partial computations are possible.
\end{rmk}


The key fibration $p^{A_1}_{\mu_1}$ is
$$X_{\mu_1}\xra{\Gamma\CP^\infty}\Gamma T\tilde\xi_{A_1}.$$
By \cite{rim} the vector bundle $\tilde\xi_{A_1}$ has base space $\CP^\infty\times\CP^\infty$ and is of the form $\pr_1^*\gamma_1^{\otimes2}\oplus\pr_2^*\gamma_1\oplus\pr_1^*\gamma_1^\vee\otimes\pr_2^*\gamma_1$ where $\pr_1$ and $\pr_2$ denote the projections $\CP^\infty\times\CP^\infty\to\CP^\infty$ to the first and second factor respectively and $\gamma_1^\vee$ is the dual bundle of $\gamma_1$. In a more comfortable form this vector bundle is $(\gamma_1^{\otimes2}\times\gamma_1)\oplus(\pr_1^*\gamma_1^\vee\otimes\pr_2^*\gamma_1)$.

\begin{rmk}
  \label{rmk:ts+}
  For any two vector bundles $\alpha$ and $\beta$ the Thom space $T(\alpha\times\beta)$ is $T\alpha\wedge T\beta$. If the base spaces of $\alpha$ and $\beta$ are the same, then the Thom space $T(\alpha\oplus\beta)$ is $T(\alpha|_{D\beta})/T(\alpha|_{S\beta})$ where, for simplicity of notation, we did not indicate the pullback of the vector bundle $\alpha$ to the disk and sphere bundles of $\beta$.
\end{rmk}

Our goal will now be to describe the Thom space appearing in the key fibration, i.e. the space
$$T\tilde\xi_{A_1}\cong T((\gamma_1^{\otimes2}\times\gamma_1)|_{D(\pr_1^*\gamma_1^\vee\otimes\pr_2^*\gamma_1)})/T((\gamma_1^{\otimes2}\times\gamma_1)|_{S(\pr_1^*\gamma_1^\vee\otimes\pr_2^*\gamma_1)}).$$

\begin{prop}
  \label{prop:txi}
  The Thom space $T\tilde\xi_{A_1}$ is homotopy equivalent to $\CP^\infty/\RP^\infty\wedge\CP^\infty/\CP^1$.
\end{prop}

To prove this proposition we will need the following three lemmas.


\begin{lemma}
  \label{lemma:tensorsph}
  For any $n\ne0$, the total space of the sphere bundle $S\gamma_1^{\otimes n}$ is $B\Z_{|n|}$.\footnote{As usual, we put $\gamma_1^{\otimes n}:=(\gamma_1^\vee)^{\otimes-n}$ if $n$ is negative.} This implies that the boundary homomorphism $\partial^n\colon\pi_2(\CP^\infty)\cong\Z\to\Z\cong\pi_1(S^1)$ in the homotopy long exact sequence of the fibration $S\gamma_1^{\otimes n}\xra{S^1}\CP^\infty$ is the multiplication by $n$ or $-n$; if we declare $\partial^1$ to be $1$, then $\partial^n$ will be $n$ for any $n\ne0$.  
\end{lemma}

\begin{prf}
  First assume that $n$ is positive. The tensor power $\gamma_1^{\otimes n}$ is obtained from
  $$\gamma_1=B\U(1)\utimes{\U(1)}\C$$
  by changing the action of its structure group $\U(1)\cong S^1$ (seen as the unit complex numbers) on $\C$ from multiplication by $\varepsilon\in S^1$ to multiplication by $\varepsilon^n$. Hence the total space of $\gamma_1^{\otimes n}$ is the factor of the total space of $\gamma_1$ by the $\Z_n$-action which fibrewise multiplies vectors by the $n$'th roots of unity. The restriction of this quotient map to the sphere bundle $S\gamma_1\cong S^\infty\cong*$ is a principal $\Z_n$-bundle $S\gamma_1\xra{\Z_n}S\gamma_1^{\otimes n}$ with contractible total space, hence $S\gamma_1^{\otimes n}$ can only be $B\Z_n$. Now if $n$ was negative, then this reasoning can be repeated with $\gamma_1^\vee$ replacing $\gamma_1$.

  This also implies that for any $n>0$ we have a commutative diagram of fibrations
  $$\xymatrix@R=1pc{
    S\gamma_1^{\otimes n}\ar[dd]^{S^1} & S\gamma_1\ar[dd]^{S^1}\ar[l]\ar[r] & S\gamma_1^\vee\ar[dd]^{S^1}\ar[r] & S\gamma_1^{\otimes-n}\ar[dd]^{S^1}\\
    & \ar[l]+<1.333pc,0pc>_n & \ar@{<-}[l]+<1.333pc,0pc>_{-1} & \ar@{<-}[l]+<1.333pc,0pc>_n\\
    \CP^\infty & \CP^\infty\ar[l]_\id\ar[r]^\id & \CP^\infty\ar[r]^\id & \CP^\infty
  }$$
  where the numbers on the arrows between the fibres $S^1$ mean maps inducing on $\pi_1(S^1)\cong\Z$ multiplications by them. Taking homotopy groups here yields a commutative diagram connecting the long exact sequences corresponding to these fibrations. Thus the boundary homomorphisms fit into a commutative diagram
  $$\xymatrix{
    \pi_2(\CP^\infty)\ar[d]_{\partial^n} & \pi_2(\CP^\infty)\ar[d]_{\partial^1}^1\ar[l]_\id\ar[r]^\id & \pi_2(\CP^\infty)\ar[d]_{\partial^{-1}}\ar[r]^\id & \pi_2(\CP^\infty)\ar[d]_{\partial^{-n}}\\
    \pi_1(S^1) & \pi_1(S^1)\ar[l]_n\ar[r]^{-1} & \pi_1(S^1) \ar[r]^n & \pi_1(S^1)
  }$$
  which concludes our proof.
\end{prf}

\begin{lemma}
  \label{lemma:diag}
  The total space of $S(\pr_1^*\gamma_1^\vee\otimes\pr_2^*\gamma_1)$ is homotopy equivalent to $\CP^\infty$ and its projection to $\CP^\infty\times\CP^\infty$ is homotopic to the diagonal map $x\mapsto(x,x)$.
\end{lemma}

\begin{prf}
  Putting $S:=S(\pr_1^*\gamma_1^\vee\otimes\pr_2^*\gamma_1)$, the non-trivial segment of the homotopy long exact sequence of the fibration $p\colon S\xra{S^1}\CP^\infty\times\CP^\infty$ is
  $$0\to\pi_2(S)\to\pi_2(\CP^\infty\times\CP^\infty)\to\pi_1(S^1)\to\pi_1(S)\to0.$$
  Here we have $\pi_2(\CP^\infty\times\CP^\infty)\cong\pi_2(\CP^\infty\times*)\oplus\pi_2(*\times\CP^\infty)\cong\Z\oplus\Z$ (for a point $*\in\CP^\infty$) where the isomorphism is induced by the inclusions $\CP^\infty\times*\subset\CP^\infty\times\CP^\infty\supset*\times\CP^\infty$. The restrictions of $S$ over $\CP^\infty\times*$ and $*\times\CP^\infty$ are $S\gamma_1^\vee$ and $S\gamma_1$ respectively, hence using lemma \ref{lemma:tensorsph} we get that the boundary homomorphism $\pi_2(\CP^\infty\times\CP^\infty)\to\pi_1(S^1)$ is generated by $(1,0)\mapsto-1$, $(0,1)\mapsto1$.

  Thus $\pi_1(S)$ is trivial and $\pi_2(S)\cong\Z$ yielding that $S$ is $K(\Z,2)\cong\CP^\infty$. The image of $p_*\colon\pi_2(S)\cong\Z\to\Z\oplus\Z\cong\pi_2(\CP^\infty\times\CP^\infty)$ is generated by $(1,1)$, hence we can identify $\pi_2(S)$ with $\Z$ such that $p_*$ is given by $1\mapsto(1,1)$. Now since $\CP^\infty\times\CP^\infty$ is $K(\Z\oplus\Z,2)$, the map $p\colon\CP^\infty\to\CP^\infty\times\CP^\infty$ is determined up to homotopy by its induced map in $\pi_2$ meaning that it is indeed homotopic to the diagonal map.
\end{prf}

Now observe that if $\Delta\subset\CP^\infty\times\CP^\infty$ denotes the diagonal $\{(x,x)\mid x\in\CP^\infty\}$, then the restriction of the projection $\pr_1$ to the subspace $E:=\CP^\infty\times\CP^\infty\setminus\Delta$ is a fibration over $\CP^\infty$ with fibre $\CP^\infty\setminus*$ (where $*$ is a point in $\CP^\infty$). More precisely the fibre $E_p$ over a point $p\in\CP^\infty$ is $\CP^\infty\setminus\{p\}$; in this fibre let $B_p$ be the subspace $\{\ell<\C^\infty\mid\dim\ell=1,\ell\perp p\}$ and let $B\subset\CP^\infty\times\CP^\infty\setminus\Delta$ be the union of the spaces $B_p$ for all points $p\in\CP^\infty$. For $i=1,2$ denote by $\tilde\pr_i$ the restriction $\pr_i|_B$. Note that $E$ is the total space of a line bundle over $B$; denote this bundle by $q\colon E\to B$

\begin{lemma}
  \label{lemma:minusdiag}
  The space $B$ is homotopy equivalent to $\CP^\infty\times\CP^\infty$ and the maps $\tilde\pr_i\colon B\to\CP^\infty$ are homotopic to the original projections $\pr_i$ (for $i=1,2$), moreover, the line bundle $q$ is isomorphic to $\pr_2^*\gamma_1$ (which, of course, is the same as $\pr_1^*\gamma_1$).
\end{lemma}

\begin{prf}
  The map $\tilde\pr_1\colon B\to\CP^\infty$ is a fibration with fibre $\CP^{\infty-1}$; the non-trivial segment of its homotopy long exact sequence is
  $$0\to\pi_2(\CP^{\infty-1})\to\pi_2(B)\to\pi_2(\CP^\infty)\to0.$$
  Thus $B$ can only be $K(\Z\oplus\Z,2)\cong\CP^\infty\times\CP^\infty$ and the map $\tilde\pr_1$ is homotopic to the projection $\pr_1$. The inclusion of a fibre $B_p\cong\CP^{\infty-1}\subset B$ is homotopic to the inclusion of a fibre of the trivial fibration $\pr_1$. Since $\CP^{\infty-1}$ is a deformation retract of $\CP^\infty$, the other projection $\pr_2\colon B\to\CP^\infty$ of this trivial fibration is also homotopic to the map $\tilde\pr_2$.

  For any fibre $B_p$ of $\tilde\pr_1$ the restriction of the line bundle $q\colon E\to B$ over $B_p\cong\CP^{\infty-1}$ is the pullback $\tilde\pr_2|_{B_p}^*\gamma_1$. Hence we also get that $q$ is isomorphic to $\tilde\pr_2^*\gamma_1\cong\pr_2^*\gamma_1$.
\end{prf}


\begin{prf}[of proposition \ref{prop:txi}]
  By remark \ref{rmk:ts+} and lemma \ref{lemma:diag} we have
  \begin{alignat*}2
    T\tilde\xi_{A_1}&\cong T((\gamma_1^{\otimes2}\times\gamma_1)|_{D(\pr_1^*\gamma_1^\vee\otimes\pr_2^*\gamma_1)})/T((\gamma_1^{\otimes2}\times\gamma_1)|_{S(\pr_1^*\gamma_1^\vee\otimes\pr_2^*\gamma_1)})\cong\\
    &\cong T(\gamma_1^{\otimes2}\times\gamma_1)/T((\gamma_1^{\otimes2}\times\gamma_1)|_{\Delta}).
  \end{alignat*}
  Lemma \ref{lemma:minusdiag} implies that the space $T(\gamma_1^{\otimes2}\times\gamma_1)\setminus T((\gamma_1^{\otimes2}\times\gamma_1)|_{\Delta})$ coincides up to homotopy with the Thom space of $\gamma_1^{\otimes2}\times\gamma_1$ restricted to the open disk bundle of $\pr_2^*\gamma_1$, moreover, taking one-point compactifications yields that the factor of $T(\gamma_1^{\otimes2}\times\gamma_1)$ with the subspace $T((\gamma_1^{\otimes2}\times\gamma_1)|_{\Delta})$ is the same as the factor of $T((\gamma_1^{\otimes2}\times\gamma_1)|_{D(\pr_2^*\gamma_1)})$ with $T((\gamma_1^{\otimes2}\times\gamma_1)|_{S(\pr_2^*\gamma_1)})$. Hence the above factor space is
  $$T\tilde\xi_{A_1}\cong T((\gamma_1^{\otimes2}\times\gamma_1)|_{D(\pr_2^*\gamma_1)})/T((\gamma_1^{\otimes2}\times\gamma_1)|_{S(\pr_2^*\gamma_1)})$$
  which is the Thom space of $(\gamma_1^{\otimes2}\times\gamma_1)\oplus\pr_2^*\gamma_1=
  \pr_1^*\gamma_1^{\otimes2}\oplus2\pr_2^*\gamma_1=\gamma_1^{\otimes2}\times2\gamma_1$. Thus we have
  $$T\tilde\xi_{A_1}\cong T(\gamma_1^{\otimes2}\times2\gamma_1)\cong T\gamma_1^{\otimes2}\wedge T(2\gamma_1)\cong\CP^\infty/\RP^\infty\wedge\CP^\infty/\CP^1$$
  using remark \ref{rmk:ts+} and lemma \ref{lemma:tensorsph} in the second and third equivalences.
\end{prf}

Thus we obtained the following:

\begin{thm}
  \label{thm:foldkey}
  The key fibration $p_{\mu_1}^{A_1}$ is of the form
  $$X_{\mu_1}\xra{\Gamma\CP^\infty}\Gamma(\CP^\infty/\RP^\infty\wedge\CP^\infty/\CP^1),$$
  hence there is a long exact sequence
  $$\ldots\to\pi^s_{n+2}(\CP^\infty)\to\Qhol_{\mu_1}(n,1)\to\pi^s_{n+2}(\CP^\infty/\RP^\infty\wedge\CP^\infty/\CP^1)\to\pi^s_{n+1}(\CP^\infty)\to\ldots$$
\end{thm}

As an application of this we 
can determine (the vanishing of) some of the torsion parts of the 
cobordism groups of quasi-holomorphic fold maps $M^n\to\R^{n+2}$ if $n$ is relatively small. 

\begin{crly}
  For any prime $p$ and any $n<4p-5$ the $p$-primary torsion of the cobordism group $\Qhol_{\mu_1}(n,1)$ is isomorphic to that of $\pi^s_{n+2}(\CP^\infty/\RP^\infty\wedge\CP^\infty/\CP^1)$.
\end{crly}

\begin{prf}
  We saw in theorem \ref{thm:keyractrivi} that the key fibration is rationally trivial, meaning that the above long exact sequence rationally splits to short exact sequences with $\Qhol_{\mu_1}(n,1)\otimes\Q$ as the middle term. This implies that the torsion elements in $\Qhol_{\mu_1}(n,1)$ either come from torsion elements in $\pi^s_{n+2}(\CP^\infty)$ or go to torsion elements in $\pi^s_{n+2}(\CP^\infty/\RP^\infty\wedge\CP^\infty/\CP^1)$. It was shown in \cite{liu} that for any prime $p$ the group $\pi^s_{n+2}(\CP^\infty)$ has no $p$-primary torsion if $n+2<4p-3$, hence in this range of dimensions the $p$-primary torsion part of $\Qhol_{\mu_1}(n,1)$ is isomorphic to that of $\pi^s_{n+2}(\CP^\infty/\RP^\infty\wedge\CP^\infty/\CP^1)$.
\end{prf}

\begin{crly}
  \label{crly:foldtors}
  For any $n$ the cobordism group $\Qhol_{\mu_1}(n,1)$ has no $p$-primary torsion if $p$ is a prime greater than $\frac{n+5}2$.
\end{crly}

\begin{prf}
  We first observe that the prime $p=2$ is automatically excluded, hence we may assume $p>2$ from now on. Let $[-]_p$ denote the $p$-localisation functor. Since all homology groups of the space $\RP^\infty$ are finite $2$-primary, we have $[\RP^\infty]_p\cong*$, hence we also have $[\CP^\infty/\RP^\infty\wedge\CP^\infty/\CP^1]_p\cong[\CP^\infty\wedge\CP^\infty/\CP^1]_p$.

  The Atiyah--Hirzebruch spectral sequence for the stable homotopy groups of $[\CP^\infty\wedge\CP^\infty/\CP^1]_p$ has starting page $E^2_{i,j}=H_i([\CP^\infty\wedge\CP^\infty/\CP^1]_p;\pi^s(j))$ and converges to the $p$-localised stable homotopy groups $[\pi^s_{i+j}(\CP^\infty\wedge\CP^\infty/\CP^1)]_p$ which, for $i+j<4p-3$ are the same as $[\Qhol_{\mu_1}(i+j-2,1)]_p$ by the previous corollary. By a theorem of Serre \cite{serre} the $p$-primary part of the stable homotopy group of spheres $\pi^s(j)$ is trivial for $0<j<2p-3$ and it is $\Z$ for $j=0$ and $\Z_p$ for $j=2p-3$. Thus for $i+j<2p-3$ we have $E^2_{i,j}\cong E^\infty_{i,j}$ which is $[\Z]_p^{\frac i2-2}$ if $j=0$ and $i$ is even and trivial otherwise. Hence for $n<2p-5$ the groups $[\pi^s_{n+2}(\CP^\infty\wedge\CP^\infty/\CP^1)]_p\cong[\Qhol_{\mu_1}(n,1)]_p$ are trivial.
\end{prf}






\begin{ex}
  For the first few dimensions $n$ we can completely determine the cobordism group of quasi-holomorphic fold maps $M^n\to\R^{n+2}$:
  \begin{table}[H]
    \begin{center}
      \begin{tabular}{c||c|c|c|c|c|c}
        $n$ & $0$ & $1$ & $2$ & $3$ & $4$ & $5$ \\
        \hline
        $\Qhol_{\mu_1}(n,1)$ & $\Z$ & $0$ & $\Z$ & $0$ & $\Z\oplus\Z$ & $\Z_2$ or $0$.
      \end{tabular}
    \end{center}
  \end{table}
  \vspace{-.5cm}\noindent
  Since the (complex) codimension of $A_1$ in these dimensions is $2$, for $n=0,1,2$ neither the quasi-holomorphic maps $M^n\to\R^{n+2}$, nor the cobordisms $W^{n+1}\to\R^{n+2}\times[0,1]$ between them contain singularities, hence then we have $\Qhol_{\mu_1}(n,1)\cong\Imm^{\gamma_1}(n,2)\cong\pi^s_{n+2}(\CP^\infty)$. For larger dimensions the homotopy long exact sequence of the key fibration is as follows:
  \begin{alignat*}2
    \ldots\to\pi^s_7(\CP^\infty)\cong\Z_2\to\Qhol_{\mu_1}(5,1)&\to\pi^s_7(T\tilde\xi_{A_1})\to\pi^s_6(\CP^\infty)\cong\Z\to\Qhol_{\mu_1}(4,1)\to\\
          &\to\pi^s_6(T\tilde\xi_{A_1})\to\pi^s_5(\CP^\infty)\cong\Z_2\to\Qhol_{\mu_1}(3,1)\to0
  \end{alignat*}
  where we know the stable homotopy groups of $\CP^\infty$ from \cite{liu}. To get the above stated results on $\Qhol_{\mu_1}(n,1)$, it is enough to show that we have $\pi^s_7(T\tilde\xi_{A_1})\cong0$ and $\pi^s_6(T\tilde\xi_{A_1})\cong\Z$ and the boundary map $\partial\colon\pi^s_6(T\tilde\xi_{A_1})\to\pi^s_5(\CP^\infty)$ in the sequence is non-trivial.

  From \cite[lemma 1]{2k+2} we immediately get $\pi^s_6(T\tilde\xi_{A_1})\cong\Z$ since $\tilde\xi_{A_1}$ is an orientable $6$-vector bundle. Moreover, since the base space $\CP^\infty\times\CP^\infty$ of $\tilde\xi_{A_1}$ is simply connected and we have
  \begin{alignat*}2
    w_2(\tilde\xi_{A_1})&=w_2(\textstyle\pr_1^*\gamma_1^{\otimes2}\oplus\pr_2^*\gamma_1\oplus\pr_1^*\gamma_1^\vee\otimes\pr_2^*\gamma_1)=\\
    &=w_2(\textstyle\pr_1^*\gamma_1^{\otimes2})+w_2(\pr_2^*\gamma_1)+w_2(\pr_1^*\gamma_1^\vee)+w_2(\pr_2^*\gamma_1)=w_2(\pr_1^*\gamma_1^\vee)\ne0,
  \end{alignat*}
  we also get $\pi^s_7(T\tilde\xi_{A_1})\cong0$ from \cite[lemma 2]{2k+2}.

  Now observe that the group $\pi^s_6(T\tilde\xi_{A_1})\cong\Imm^{\tilde\xi_{A_1}}(0,6)\cong\Z$ is generated by the cobordism class of a point embedded into $\R^6$ with normal bundle induced from $\tilde\xi_{A_1}$. By its construction, the boundary map $\partial$ maps this generator to the cobordism class in $\pi^s_5(\CP^\infty)\cong\Imm^{\gamma_1}(3,2)\cong\Z_2$ of the immersion $f\colon S^3\imto\R^5$ that we get by restricting the prototype of the fold singularity $\C^2\to\C^3$ to the preimage of a sufficiently small sphere (say, of radius $\varepsilon$) around the origin. This is a map of a manifold diffeomorphic to $S^3$ to $S^5_\varepsilon$ (which is up to cobordism the same as a map to $\R^5$) and it is called the \textit{link} of $A_1$. Now the results in \cite{npimm} show that the link of $A_1$ is not null-cobordant as an immersion, since 
  otherwise the \textit{cross-cap number} $C(A_1)$ would be $0$ but it is known to be $1$. Thus we obtained that the boundary homomorphism $\partial$ is non-trivial. This is what we wanted.

  Lastly we note that to decide whether the group $\Qhol_{\mu_1}(5,1)$ is $0$ or $\Z_2$ we should know whether an immersion whose cobordism class generates the group $\pi^s_7(\CP^\infty)\cong\Imm^{\gamma_1}(5,2)\cong\Z_2$ is null-cobordant as a quasi-holomorphic map or not. That is, if $f\colon M^5\imto\R^7$ represents such a generator, is there a manifold $W^6$ with boundary $\partial W=M$ and a quasi-holomorphic map $F\colon W\to\R^8_+$ which restricts to $f$ on the boundary (for dimensional reasons, if such a map exists, it is automatically a fold map).
\end{ex}

\section{Further directions}
\label{sec:fin}

In this final section we gather a few directions in which we may move forward with the study of quasi-holomorphic maps and their applications. We note here that some of these topics 
arose from discussions with Gergő Pintér.
\begin{enumerate}[label=\textbf{\Roman*.}]
\item\textit{Links of finitely determined germs:} Let $\eta\colon(\C^n,0)\to(\C^{n+k},0)$ be a finitely determined holomorphic germ which is a prototype of its singularity class. Then the link of $\eta$, defined as the restriction of $\eta$ to the preimage of a sufficiently small sphere $S^{2(n+k)-1}_\varepsilon$ around the origin, is a stable quasi-holomorphic map $f:=\eta|_{\eta^{-1}(S^{2(n+k)-1}_\varepsilon)}\colon S^{2n-1}\to S^{2(n+k)-1}$ (both the source and the target spaces can be identified with spheres). The link of a germ is a widely studied, important invariant in local singularity theory; in fact it contains all topological information on $\eta$ (see e.g. the textbook \cite{mondnb} for general references). Now if $\sigma$ is the set of singularities $\vartheta\le\eta$ and we put $\sigma':=\sigma\setminus\{\eta\}$, then global singularity theoretic properties of $f$ as a quasi-holomorphic $\sigma'$-map will be topological invariants of the singularity $[\eta]$.

  As an example, take the cobordism class of $f$ which is an element $[f]$ in $\Qhol_{\sigma'}(2n-1,k)$. The homotopy long exact sequence of the key fibration $p_\sigma^\eta$ includes the portion
  $$
  \pi^s_{2(n+k)}(T\tilde\xi_\eta)\xra\alpha\Qhol_{\sigma'}(2n-1,k)\xra\beta\Qhol_\sigma(2n-1,k)$$
  and since $\eta$ itself defines a null-cobordism of $f$ as a $\sigma$-map, we have $[f]\in\ker\beta=\im\alpha$. Actually it is not hard to see that $\pi^s_{2(n+k)}(T\tilde\xi_\eta)\cong\Z$ (e.g. see \cite[lemma 1]{2k+2}) and $[f]$ is the image of the generator of this $\Z$, moreover, corollary \ref{crly:rachol} implies that $\beta$ is a rational isomorphism, hence $[f]$ is a torsion element in $\Qhol_{\sigma'}(2n-1,k)$. Thus for example the order of $[f]$ gives information on the topology of $\eta$.
\item\textit{Immersion formulae:} In the early 2000's, starting with works of Ekholm and the second author (see e.g. \cite{ekholm}, \cite{esz}, \cite{sszt}), se\-veral invariants of regular homotopy classes of immersions $M\imto P$ (with various additional conditions; usually with $M$ being a sphere and $P$ a Eucledian space) were constructed, essentially by considering its multiple points and also the multisingularities of generic mappings of a null-cobordism of the manifold $M$. Now if we consider an immersion $f\colon M^n\imto P^{n+2k}$ with a complex structure on its normal bundle, where $M$ is a null-cobordant manifold, then we can ask whether there is a quasi-holomorphic map from a null-cobordism of $M$ extending $f$. If there is, then the multisingularities of such a map may also yield invariants of the regular homotopy class of $f$.

  One of the simplest examples where singularities occur is that of an immersion $f\colon S^3\imto\R^5$ whose normal bundle we endow with a complex structure (this can always be done since for a rank-$2$ vector bundle a complex structure is the same as an orientation). If $W^4$ is a manifold with boundary $\partial W=S^3$, then the only possible multisingularities of a stable quasi-holomorphic map $W\to\R^6_+$ extending $f$ are $1\Sigma^0$ (on a $4$-dimensional subspace), $2\Sigma^0$ (on a $2$-dimensional subspace), $3\Sigma^0$ (in discrete points) and $1A_1$ (in discrete points). It may turn out that (similarly to the real case) some combination of invariants of these multisingular loci is an invariant that only depends on the immersion $f$, thus we can use it to gain information on $f$. We note that this is connected to the previous topic since $f$ can appear as the link of a finitely determined germ $\C^2\to\C^3$, as in e.g. \cite{npimm}.
\item\textit{Representing (co)homology classes:} Thom \cite{th} proved that for any space $P$, every homology class $x\in H_n(P;\Z_2)$ can be represented by the mapping of a manifold $f\colon M^n\to P$ (i.e. we have $x=f_*[M]$). If $P$ is a manifold of dimension $n+2k$, then by Poincaré duality, this is equivalent to saying $f_!1_M=y\in H^{2k}(P;\Z_2)$ where $1_M$ is the fundamental class of $M$ and $y$ is the Poincaré dual of $x$. Now one can ask for conditions under which the map $f$ representing the cohomology class $y$ can be chosen to be quasi-holomorphic. That is, for a manifold $P^{n+2k}$, which cohomology classes $y\in H^{2k}(P;\Z_2)$ are such that there is a quasi-holomorphic map $f\colon M^n\to P$ with $f_!1_M=y$? Note that since quasi-holomorphic maps are naturally equipped with normal orientations (a corollary of lemma \ref{lemma:cplxnb}), we know from \cite[theorem B]{grantsz} that this cannot hold for all $y$. Also note that by the same normal orientation the cohomology class $f_!1_M$ might also exist with integer (or rational, etc.) coefficients, which motivates the following: For a manifold $P^{n+2k}$ and a coefficient group $G$, which cohomology classes in $H^{2k}(P;G)$ 
  can be represented by a quasi-holomorphic map in the same sense as above?

  A more restrictive connection of the complex world to the problem of representing cohomology classes by ``nice'' subspaces is the following: For a complex projective manifold $P^{n+k}$ (i.e. a smooth subvariety in some $\CP^L$), which cohomology classes $y\in H^{2k}(P;G)$ can be written as $G$-linear combinations of classes represented by complex subvarieties in $P$? Note that since subvarieties can be locally parametrised by holomorphic maps, all classes $y$ for which the answer is positive are also represented by quasi-holomorphic maps. The famous Hodge conjecture (e.g. see \cite{hodgeconj}) asserts that a cohomology class $y\in H^{2k}(P;\Q)$ is a rational linear combination of classes represented by subvarieties if and only if $y$, as a de Rham cohomology class, is represented as a harmonic form of type $(k,k)$, i.e. if $y$ is a \textit{Hodge class}. The ``only if'' part of this statement is not hard to see, but the ``if'' part is a very hard open problem. Now a simplified version of this problem is to show that for a complex projective manifold $P$, all Hodge classes of $P$ can be represented by quasi-holomorphic maps.
\item\textit{Cobordisms of branched coverings:} Branched coverings are special cases of quasi-holomorphic maps of codimension $0$, and their cobordism groups are interesting objects in themselves. Since the proof of theorem \ref{thm:key} does not work for codimension-$0$ maps, the main computational tool we used in section \ref{sec:cob} does not exist for branched coverings. However, the construction of the classifying space $X_\tau$, where $\tau$ is a set of multisingularities of branched coverings, is relatively simple (see also \cite{nagy} where it was constructed independently), hence there might be hope to compute such cobordism groups as well with some different technique. Related to this, we ask here whether theorem \ref{thm:key} can be extended to the case $k=0$ (cf. remark \ref{rmk:key}). Note that the real smooth analogue of this question is also open.
\item\textit{Complex analogue of Kazarian's conjecture:} One of the main results in the paper \cite{hosszu} was the proof of a conjecture of Kazarian identifying, up to homotopy equivalence, the real smooth version of the space $X_\sigma$ (for a set of singularities $\sigma$) with the infinite loop space $\Gamma T\nu_\sigma$ of a certain Thom spectrum given by a virtual vector bundle over the real smooth version of the Kazarian space $K_\sigma$. The holomorphic analogue of the space $\Gamma T\nu_\sigma$ can be defined in a completely analogous way to the real smooth one, and it seems likely that it is also homotopy equivalent to the holomorphic version of $X_\sigma$. The proof in \cite{hosszu} is quite involved and consists of many steps that should all be checked in the quasi-holomorphic setting as well, but we suspect that it goes through for quasi-holomorphic maps without much additional work.
\end{enumerate}
A few additional, naturally arising questions are the following:
\begin{enumerate}[label=\textbf{\Roman*.}]
\setcounter{enumi}{5}
\item \textit{For which pairs of manifolds $M^n$, $P^{n+2k}$ does a quasi-holomorphic map $M\to P$ exist?}
\item \textit{For a smooth map $f$ between smooth manifolds, when can we endow $f$ with a quasi-holomorphic structure?}
\item \textit{For a quasi-holomorphic map $f$ between complex analytic manifolds, when is $f$ quasi-holo\-morphically homotopic/(pseudo)isotopic/cobordant to a holomorphic map?}
\item \textit{When are two holomorphic maps quasi-holomorphically homotopic/(pseudo)isotopic/cobord\-ant to each other?}
\end{enumerate}


\begin{thebibliography}{0000000}{\footnotesize

  \bibitem[An08]{ando} Y. Ando, {\it Cobordisms of maps with singularities of given class}, Algebr. Geom. Topol. 8 (2008), 1989--2029.
    
  \bibitem[Bo67]{boar}  J. M. Boardman, {\it Singularities of differentiable maps}, Publ. Math. I.H.É.S. 33 (1967), 21--57.
    
  \bibitem[BK72]{bk} A. K. Bousfield, D. M. Kan, {\it Homotopy limits, completions and localizations}, Lecture Notes in Math. 304, Springer (1972).
    
  \bibitem[Ek01]{ekholm} T. Ekholm, \textit{Differential 3-knots in 5-space with and without self-intersections}, Topology 40 (2001), 157--196.
    
  \bibitem[ESz03]{esz} T. Ekholm, A. Szűcs, \textit{Geometric formulas for Smale invariants of codimension two immersions}, Topology 42 (2003), 171--196.
    
  \bibitem[GSz15]{grantsz} M. Grant, A. Szűcs, \textit{Homologies are infinitely complex}, Topol. Methods Nonlinear Anal. 45 (2015), 55--61.
    
  \bibitem[HK56]{thomthm} A. Haefliger, A. Kosinski, {\it Un théorème de Thom sur les singularités des applications différentiables}, Séminaire Henri Cartan 9 (1956/57), t. 8.
    
  \bibitem[Ka01]{kazspace} M. É. Kazarian, {\it Classifying spaces of singularities and Thom polynomials}, in: New Developments in Singularity Theory (Cambridge, 2000), NAII 21, Springer (2001), 117--134.
    
  \bibitem[Ka06]{kaztp} M. É. Kazarian, {\it Thom polynomials}, in: Singularity theory and its applications, Adv. Stud. Pure Math. 43 (2006), 85--135.

  \bibitem[Ka]{morchar} M. É. Kazarian, \textit{Morin maps and their characteristic classes}, unpublished preprint.
    
  \bibitem[Le99]{hodgeconj} J. D. Lewis, \textit{A survey of the Hodge conjecture}, CRM Monogr. Ser. 10, Amer. Math. Soc. (1999).
    
  \bibitem[Li63]{liu} A. Liulevicius, {\it A theorem in homological algebra and stable homotopy of projective spaces}, Trans. Amer. Math. Soc. 109 (1963), 540--552.

  \bibitem[Ma70]{math5} J. N. Mather, \textit{Stability of $C^\infty$ mappings, V: Transversality}, Advances in Math. 4 (1970), 301--336.
    
    
    
  \bibitem[MNB20]{mondnb} D. Mond, J. J. Nuño-Ballesteros, \textit{Singularities of mappings—the local behaviour of smooth and complex analytic mappings}, Grundlehren Math. Wiss. 357, Springer (2020).

  \bibitem[Mor65]{mor} B. Morin, {\it Formes canoniques des singularités d'une application différentiable}, C. R. Acad. Sci. Paris 260 (1965), 5662--5665 and 6503--6506.
    
  \bibitem[Mo68]{mosh} R. E. Mosher, {\it Some stable homotopy of complex projective space}, Topology 7 (1968), 179--193.

  \bibitem[Mu93]{muk} J. Mukai, {\it On stable homotopy of the complex projective space}, Japan J. Math 19 (1993), 191--216.
    
  \bibitem[Na17]{nagy} Cs. Nagy, \textit{Cobordism groups of simple branched coverings}, Acta Math. Hungar. 153 (2017), no. 2, 449--489.
    
  \bibitem[NP15]{npimm} A. Némethi, G. Pintér, \textit{Immersions associated with holomorphic germs}, Comment. Math. Helv. 90 (2015), 513--541.
    
  \bibitem[Ri01]{rim} R. Rimányi, {\it Thom polynomials, symmetries and incidences of singularities}, Invent. Math. 143 (2001), 499--521.
    
  \bibitem[Ri02]{rlsym} R. Rimányi, \textit{On right-left symmetries of stable singularities}, Math. Z. 242 (2002), 347--366.

  \bibitem[RSz98]{rsz} R. Rimányi, A. Szűcs, {\it Pontrjagin--Thom-type construction for maps with singularities}, Topology 37 (1998), 1177--1191.

  \bibitem[RS01]{rs} C. Rourke, B. Sanderson, {\it The compression theorem I}, Geom. Topol. 5 (2001), 399--429.

  \bibitem[Sa09]{sad} R. Sadykov, {\it Bordism groups of solutions to differential relations}, Algebr. Geom. Topol. 9 (2009), 2311--2347.
    
  \bibitem[SSzT02]{sszt} O. Saeki, A. Szűcs, M. Takase, \textit{Regular homotopy classes of immersions of 3-manifolds into 5-space}, Manuscripta Math. 108 (2002), 13--32.
    
  \bibitem[SO18]{sasohm} T. Sasajima, T. Ohmoto, \textit{Thom polynomials in $\AA$-classification I: counting singular projections of a surface}, in: Schubert varieties, equivariant cohomology and characteristic classes--IMPANGA 15, EMS Ser. Congr. Rep. (2018), 237--259.
    
  \bibitem[Se51]{serre} J.-P. Serre, \textit{Homologie singulière des espaces fibrés}, Ann. of Math. 54 (1951), 425--505.

  \bibitem[Sz76]{limm} A. Szűcs, \cyrins{\textit{Группы кобордизмов $l$-погружений. I}} (Cobordism groups of $l$-immersions. I), Acta Math. Acad. Sci. Hungar. 27 (1976), 343--358.

  \bibitem[Sz93]{univmap} A. Szűcs, \textit{Universal singular map}, in: Topology. Theory and applications, II (Pécs, 1989), Colloq. Math. Soc. János Bolyai, 55 (1993), 491--500.
    
  \bibitem[Sz08]{hosszu} A. Szűcs, {\it Cobordism of singular maps}, Geom. Topol. 12 (2008), 2379--2452.

  \bibitem[Te08]{2k+2} T. Terpai, {\it Cobordisms of fold maps of $2k+2$-manifolds into $\mathbb{R}^{3k+2}$}, in: Geometry and Topology of Caustics--CAUSTICS’06 (Warsaw), Banach Center Publ. 82 (2008), 209--213.

  \bibitem[Te09]{key} T. Terpai, {\it Fibration of classifying spaces in the cobordism theory of singular maps}, Proc. Steklov Inst. Math. 267 (2009), 270--277.
    
  \bibitem[Th54]{th} R. Thom, \textit{Quelques propriétés globales des variétés différentiables}, Comment. Math. Helv. 28 (1954), 17--86

  \bibitem[Th56]{tp} R. Thom, {\it Les singularités des applications différentiables}, Ann. Inst. Fourier 6 (1956), 43--87.

  \bibitem[Uch69]{eximm} F. Uchida, \textit{Exact sequences involving cobordism groups of immersions}, Osaka J. Math. 6 (1969), 397--408.
    
    
  \bibitem[Wa81]{wallfin} C. T. C. Wall, \textit{Finite determinacy of smooth map-germs}, Bull. London Math. Soc. 13 (1981), 481--539.
    
  \bibitem[We66]{wells} R. Wells, {\it Cobordism groups of immersions}, Topology 5 (1966), 281--294.

}\end{thebibliography}
\end{document}